%% file: Main.tex
\parskip=\baselineskip
\baselineskip=\baselineskip
\hoffset = 0.5 true in
\hsize 5.5  true in
\voffset = 0.5 true in
\vsize 8 true in
\def\endofproof{\hfill $\#$ } 
\def\inchap{\vskip 1.25 true in}
\def\insec{\bigskip\bigskip}
\def\tskip{\vskip -15 true pt}
\vskip -10 true pt
\font \bigfont=cmb10 at 17 true pt
\font \medfont=cmb10 at 15 true pt
\inchap
\centerline{\bigfont On Post Critically Finite Polynomials}
\centerline{\bigfont Part Two: Hubbard Trees}
\bigskip
\bigskip
\bigskip
\centerline{Alfredo Poirier}
\centerline{Math Department}
\centerline{Suny@StonyBrook}
\centerline{StonyBrook, NY 11794-3651}
\bigskip
\bigskip
\bigskip
{\sl We provide an effective classification of postcritically finite polynomials as dynamical systems by means of Hubbard Trees.} 
\bigskip
\bigskip
This paper is the second in a series of two 
based on the author's thesis 
which deals with the classification of postcritically finite polynomials as 
dynamical systems (see [P2]). 
In the first part [P1], we conclude the study of critical portraits initiated by Fisher (see [F]) 
and continued by Bielefeld, Fisher and Hubbard (see [BHF]). 
As an application of our results, 
we give in this second part of the series 
necessary and sufficient conditions for the realization of Hubbard Trees. 
\vfill\eject
\input intro.hub.tex

\vfill\eject 
\input chap.5.tex 
\vfill\eject 
\input chap.6.tex 
\vfill\eject 
\input chap7.tex

\vfill\eject 
\input apen.ht.tex
\vfill\eject
\end

%% file: intro.hub.tex
\inchap
\centerline{\bigfont Hubbard Trees.} 
\inchap
Given a polynomial $P$ of degree $n \ge 2$, 
we consider the set $K(P)$ (called the {\it filled Julia set}) of points whose orbit under iteration is bounded. 
This set is known to be compact and its complement consists of a unique unbounded component (see [M, Lemma 17.1]). 
The behavior under iteration of the critical points of this polynomial dramatically influences the topology of this set $K(P)$. 
For example, this set is connected if and only if 
all critical points are contained within (see [M, Theorem 17.3]). 
We are interested in the special case where the orbit of every critical points is finite, 
i.e, 
the case where the orbits of all critical points are periodic or eventually periodic. 
We call such polynomials {\it postcritically finite} 
($PCF$ in short). 
For such polynomials the filled Julia set $K(P)$ is connected. 
Furthermore, it is also known in this case that $K(P)$ is locally connected 
(see [M, Theorem 17.5]). 
\medskip
In order to proceed further we establish some notation. 
The set\break 
$J(P)=\partial K(P)$ is called the {\it Julia set}, and its elements {\it Julia points}. 
The complement $F(P)={\bf C}-J(P)$ of the Julia set 
is called the {\it Fatou set} and its elements {\it Fatou points}. 
A periodic orbit $z_0 \mapsto z_1=P(z_0) \mapsto \dots \mapsto z_n=z_0$ which contains a critical point is called a {\it critical cycle}. 
In the $PCF$ case a periodic orbit belongs to the Fatou set $F(P)$ if and only if it is a critical cycle (see [M, Corollary 11.6]). 
\medskip
In this $PCF$ case the dynamics of the polynomial admits a further decomposition. 
When restricted to the interior of $K(P)$ 
(which is not empty if and only if there exists a critical cycle), 
$P$ maps each component onto some other as a branched covering map. 
Furthermore, every component is eventually periodic (see [M, Theorem 13.4]). 
It is well known (see [M, Theorem 6.7]) that each component can be uniformized so that in local coordinates $P$ can be written as $z \mapsto z^n$ for some $n\ge 1$. 
Furthermore, if $U$ is a periodic bounded Fatou component, 
the first return map is conjugate in local coordinates to $z \mapsto z^n$ for some $n \ge 2$.  
In particular such cycles of components are in one to one correspondence with critical cycles. 
Also, in each component there is a unique point which eventually maps to a critical point 
(these points are those which correspond to $0$ in local coordinates).
\medskip
In the work [DH1], Douady and Hubbard suggested a combinatorial description of the dynamics of such polynomials using a tree-like structure. 
First we note the following (see [DH1, Corollary VII.4.2 p 64]). 
\bigskip
{\bf Lemma.}
{\it Let $P$ be a PCF polynomial. 
Then for any $z \in J(P)$, 
the set $J(P)-\{z\}$ consists of only a finite number of connected components.}
\bigskip
Thus, the filled Julia set is arranged in a tree like fashion. 
To simplify this tree we consider a finite invariant set $M$ (i.e, $P(M) \subset M$) containing all critical points.  
We join them in $K(P)$ by paths subject to the restriction that if they intersect a Fatou component, 
this intersection must consist of radial segments in the coordinate described above. 
Douady and Hubbard proved that this construction is unique and 
defines a finite topological tree $T_M$ in which all points in $M$ 
(and perhaps more) are vertices. 
Now, if from this tree we retain the dynamics and local degree at every vertex, 
the way this tree is embedded in the complex plane (up to isotopy class), 
and ``a bit of extra information to recover the tree generated by $P^{-1}(M)$" 
(there are several ways to state this condition in a non-ambiguous way), 
they proved that different $PCF$ polynomials 
(i.e, not conjugated as dynamical systems) 
give rise to different tree-structures. 
No criterion for realization was given at the time. 
(The only previous partial results about realization are given in 
Lavaurs' thesis [L]). 
\bigskip
A way to deal with this conditions is to introduce angles around vertices in the tree structure (see [DH1, p.46]). 
In what follows we will measure angles in turns (i.e, $360^\circ=1$ turn). 
Around a Fatou vertex $v$ 
(which correspond to $0$ in the uniformizing coordinate), 
an angle between edges incident at $v$ is naturally defined by means of the local coordinate system. 
At Julia vertices, where $m$ components of $K(P)$ meet (compare the lemma above), 
the angle is defined to be a multiple of $1/m$ 
(this normalization is introduced here for the first time). 
These angles satisfy two conditions. 
First, they are compatible with the embedding of the tree. 
Second, we have that $\angle_{P(v)} (P(\ell),P(\ell'))=\delta(v)\angle_v (\ell,\ell')$ {\it (mod 1)}, 
where $\delta(v)$ is the local degree of $P$ at $v$ and $\ell,\ell'$ are edges incident at $v$ 
($\angle_v$ and $\angle_{P(v)}$ measure the angles at $v$ and $P(v)$ respectively).
When this further structure is given, 
we have a `dynamical tree', 
which we denote by ${\bf H}_{P,M}$. 
\bigskip
Now let us start with an abstract tree and try to reconstruct the appropriate 
polynomial. 
\bigskip
{\bf Definition.} 
By an {\it (angled) tree H} will be meant a finite connected acyclic $m$-dimensional simplicial complex ($m=0,1$), 
together with a function 
$\ell,\ell' \mapsto \angle(\ell,\ell')=\angle_v(\ell,\ell') \in {\bf Q/Z}$ 
which assigns a rational modulo 1 to each pair of edges $\ell,\ell'$ which meet at a common vertex $v$. 
This angle $\angle(\ell,\ell')$ should be skew-symmetric, with $\angle(\ell,\ell')=0$ if and only if $\ell=\ell'$, 
and with $\angle_v(\ell,\ell'')=\angle_v(\ell,\ell')+\angle_v(\ell',\ell'')$ whenever three edges are incident at a vertex $v$. 
Such an angle function determines a preferred isotopy class of embeddings of $H$ into ${\bf C}$. 
\medskip
Let $V$ be the set of vertices. 
We specify a mapping $\tau:V \to V$ and call it the {\it vertex dynamics}, 
and require that $\tau(v) \ne \tau(v')$ whenever $v$ and $v'$ are endpoints of a common edge $\ell$. 
We consider also a {\it local degree function} $\delta:V \to {\bf Z}$ which assigns an integer $\delta(v) \ge 1$ to each vertex $v \in V$. 
We require that $deg(\delta)=1+\sum_{v \in V}(\delta(v)-1)$ be greater that 1.\break
By definition a vertex {\it v is critical} if $\delta(v)>1$, and 
{\it non-critical} otherwise. 
The {\it critical set} $\Omega(\delta)=\{v \in V:\hbox{\it v is critical}\}$ is thus non empty. 
\bigskip
The maps $\tau$ and $\delta$ must be related in the following way. 
Extend $\tau$ to a map $\tau:H \to H$ which carries each edge homeomorphically onto the shortest path joining the images of its endpoints. 
We require then that $\angle_{\tau(v)}(\tau(\ell),\tau(\ell'))=\delta(v)\angle_v(\ell,\ell')$
whenever $\ell,\ell'$
are incident at $v$ (in this case $\tau(\ell)$ and $\tau(\ell')$ are incident at the vertex $\tau(v)$ where the angle is measured). 
\medskip
A vertex $v$ is {\it periodic} if for some $n>0$, $\tau^{\circ n}(v)=v$. 
The orbit of a periodic critical point is a {\it critical cycle}. 
We say that a vertex $v$ is of {\it Fatou type or a Fatou vertex} if it eventually maps into a critical cycle. 
Otherwise, if it eventually maps to a non critical cycle, it is of {\it Julia type or a Julia vertex}.    
\medskip
We define the distance $d_{H}(v,v')$ between vertices in $H$ as the number of edges in a shortest path $\gamma$ between $v$ and $v'$. 
We say that $(H,V,\tau,\delta)$ is 
{\it expanding} if the following condition is satisfied. 
For any edge $\ell$ whose end points $v,v'$ are Julia vertices, 
there is an $n \ge 1$ such that $d_{H}(\tau^{\circ n}(v),\tau^{\circ n}(v'))>1$. 
\medskip
The angles at Julia vertices are rather artificial, 
so we normalize them as follows. 
If $m$ edges $\ell_1,\dots,\ell_m$ meet at a periodic Julia vertex $v$, 
then we assume that the angles $\angle_v(\ell_i,\ell_j)$ are all multiples of $1/m$. 
(It follows that the angles at a periodic Julia vertex convey no information beyond the cyclic order of these $m$ incident edges.) 
\bigskip
By an {\it abstract Hubbard Tree} we mean an angled tree 
${\bf H}=((H,V,\tau,\delta),\angle)$ 
so that the angles at any periodic Julia vertex where $m$ edges meet are multiples of $1/m$. 
\medskip
The basic existence and uniqueness theorem can now be stated as follows 
(compare Theorem II.4.7). 
\bigskip
{\bf Theorem A.}  
{\it Any abstract Hubbard Tree ${\bf H}$ can be realized as a tree associated with 
a postcritically finite 
polynomial P if and only if 
$\bf H$ is expanding. 
Such a realization is necessarily unique up to affine conjugation.} 
\bigskip
This abstract Hubbard Tree also gives information about external rays as the following theorem essentially due to Douady and Hubbard shows\break 
(compare [DH1, Chap VII]). 
This will follow in our case from Propositions II.3.3, III.4.3 and the fact
that $J(P)$ is locally connected. 
\bigskip
{\bf Theorem B.}
{\it The number of rays which land at a periodic Julia vertex $v$ is equal to 
the number of incident edges of the tree $T$ at $v$, 
and in fact, 
there is exactly one ray landing between each pair of consecutive edges. 
Furthermore, the ray which lands at $v$ between $\ell$ and $\ell'$ maps to the ray 
which lands at $f(v)$ between $f(\ell)$ and $f(\ell')$.}
\bigskip
After these theorems there is no reason to distinguish between the abstract Hubbard Tree and the unique polynomial which realizes it. 
\medskip
{\bf Definition.} 
A point $p \in J(P)$ is {\it terminal } if there is only one external ray landing at $p$. 
Otherwise $p$ is an {\it incidence point}. 
For incidence points we distinguish between {\it branching} (if there are more than two rays landing at $p$) 
and {\it non branching} (exactly two rays landing at $p$). 
For a postcritically finite polynomial $P$, 
every branching point must be periodic or preperiodic. 
Also every periodic branching point is present as vertex in any tree ${\bf H}_{P,M}$. 
\bigskip
{\bf Proposition I.3.2.} 
{\it Let P be a Postcritically Finite Polynomial and $z \in J(P)$ a 
branching point. 
Then $z$ is preperiodic (or periodic).}
\medskip
{\bf Proposition I.3.3.} 
{\it Let P be a Postcritically Finite Polynomial and $z \in J(P)$ a 
periodic incidence point. 
For any invariant finite set $M$ containing the critical points of P, 
we have $z \in T_{P,M}$. 
Furthermore, the number of components of $T_{P,M}-\{z\}$ is independent of M and equals the number of components of $J(P)-\{z\}$.}
\bigskip
\bigskip
Now we give a brief description of the following chapters. 
In Chapter I we have included the basic background of Hubbard Trees following the original exposition of 
Douady and Hubbard. 
We have done so because 
there is nowhere in the literature where we can find in a systematic way what was known up to now. 
In Chapter II, 
we introduce our basic abstract framework. 
We have carefully justified why there is the need to introduce all the 
abstract elements in our definition. 
In Chapter III we give the proof of our main result. 
This proof is based in the theory of critical portraits developed in the first part of this work. 
For the convenience of the reader we have included in Appendix A 
an outline of this theory. 
In Appendix B, 
we study necessary and sufficient conditions under which an 
$n^{th}$ fold covering of a finite cyclic set to a proper subset 
can be given a compatible `argument coordinate' so that it becomes multiplication by $n$. 
\bigskip
{\bf Acknowledgement.} 
We will like to thank John Milnor for helpful discussions and suggestions. 
Some of the arguments are in its final formulation thanks to him. 
Also, we want to thank the Geometry Center, University of Minnesota
and Universidad Cat\'olica del Per\'u 
for their material support.

%% file: chap.5.tex
\inchap
\centerline{\bigfont Chapter I}
\centerline{\bigfont Hubbard Trees}
\inchap
{In this Chapter we recall the definition and survey the main properties of Hubbard Trees 
as defined by Douady and Hubbard in [DH1]. 
In Section 1 we define the main concepts and deduce some properties. 
We ask the reader to pay special attention to Proposition 1.21. 
In Section 2 we define the inverse of Hubbard Trees. 
In Section 3 we define and study the incidence number at every point $p$ of the tree and 
relate this concept with the number of connected components of $J(P)-\{p\}$. 
} 
\insec 
\centerline{\medfont 1. Regulated Trees.}
\insec 
{\bf 1.1} 
Let $P$ be a Postcritically Finite Polynomial. 
Given two points in the closure of a bounded Fatou component, 
they can be joined in a unique way by a Jordan arc consisting of (at most two) segments of internal rays. 
We call such arcs (following Douady and Hubbard) {\it regulated}. 
The filled Julia set $K(P)$ being connected and locally connected in a compact metric space is also arcwise connected. 
This means that given two points $z_1,z_2 \in K(P)$ there is a continuous injective map $\gamma:I=[0,1]\mapsto K(P)$ such that $\gamma(0)=z_1$ and $\gamma(1)=z_2$. 
In general we will not distinguish between the map and its image.
Such arcs (actually their images) can be chosen in a unique way so that the intersection with the closure of a Fatou component is regulated
(see [DH1, Chapter 2]). 
We still call such arcs {\it regulated}, and denote them by $[z_1,z_2]_K$.
\smallskip
The following immediate properties hold for regulated arcs (compare also [DH1, Chapter 2]).
\bigskip
{\bf 1.2 Lemma.}
{\it Let $\gamma_1,\gamma_2$ be regulated arcs, then $\gamma_1 \cap \gamma_2$ is regulated.}
\endofproof
\bigskip
{\bf 1.3 Lemma.}
{\it Every subarc of a regulated arc is regulated.}
\endofproof
\bigskip
{\bf 1.4 Lemma.}
{\it Let $z_1,z_2,z_3 \in K(P)$, then there exists $p \in K(P)$ such that 
$[z_1,z_2]_K \cap [z_2,z_3]_K=[p,z_2]_K$. 
In particular if $[z_1,z_2]_K \cap [z_2,z_3]_K= \{z_2\}$, 
the set 
$[z_1,z_2]_K \cup [z_2,z_3]_K$ is a regulated arc.}
\endofproof
\bigskip
{\bf 1.5 Regulated Sets.} 
We say that a subset $X \subset K(P)$ is {\it regulated connected} 
if for every $z_1,z_2 \in X$ we have $[z_1,z_2]_K \subset X$. 
We define the {\it regulated hull $[X]_K$ of $X \subset K(P)$} as the minimal closed regulated connected subset of $K(P)$ containing $X$.
\bigskip
{\bf 1.6 Proposition.}
{\it If $z_1,...,z_n$ are points in K(P), 
the regulated hull $[z_1,...,z_n]_K$ of $\{z_1,...,z_n\}$ is a finite topological tree.}
\medskip
{\bf Proof} (Compare [DH1]). 
The proof is by induction in the number of points. 
This is clearly true for small $n$ ($=1,2$). 
Suppose $[z_1,...,z_n]_K$ is a finite topological tree, and let $z_{n+1} \in K(P)$. 
Let $p$ any point in $[z_1,...,z_n]_K$ and $y$ the first point in the arc $[z_{n+1},p]_K$ that belongs to $[z_1,...,z_n]_K$. 
In this way $[z_1,...,z_{n+1}]_K=[z_1,...,z_n]_K \cup [y,z_{n+1}]_K$ and 
$[z_1,...,z_n]_K \cap [y,z_{n+1}]_K=\{y\}$. 
The result follows. 
\endofproof
\bigskip
{\bf 1.7 Remark.} 
By definition every end of the tree $[z_1,...,z_n]_K$ is one of $z_k$, but there may be $z_k$ which are not ends.
\bigskip
{\bf 1.8 Lemma.}
{\it Let $\gamma(I) \subset K(P)$ be a regulated arc containing no critical point of P, except possibly for its end points. 
Then $P|_{\gamma(I)}$ is injective and $P(\gamma(I))$ is a regulated arc.}
\medskip
{\bf Proof.} 
The second part follows from the first, so let us show that $P|_\gamma$ is injective. 
As $P \circ \gamma$ is locally one to one, the set $\Delta=\{(t_1,t_2): t_1<t_2$ and $P(\gamma(t_1))=P(\gamma(t_2))\}$ is compact. 
If this set is non empty we can take $(t_1,t_2) \in \Delta$ with $t_2-t_1$ minimal. 
Let $t \in (t_1,t_2)$, then $P(\gamma([t_1,t]))$ and $P(\gamma([t,t_2]))$ are regulated arcs with the same end points; 
therefore they are equal and $t_2-t_1$ is not minimal.
\endofproof
\bigskip
{\bf 1.9 Definition.} 
For a finite invariant set $M$, containing 
the set $\Omega(P)$ of critical points of $P$, 
we denote by $T(M)$ the tree generated by $M$, i.e, the regulated hull $[M]_K$. 
The {\it minimal tree $T(M_0)$}, is the tree generated by $M_0={\cal O}(\Omega(P))$ the orbit of the critical set. 
This last tree is usually called in the literature {\it the Hubbard Tree of P}. 
\bigskip
{\bf 1.10 Lemma.}
{\it For a finite invariant set $M$, containing 
the set $\Omega(P)$ of all critical points, $P(T(M))=[P(M)]_K$.}
\medskip
{\bf Proof.} 
The tree $T(M)$ is the union of regulated arcs of the form $[z_1,z_2]_K$ with $z_1,z_2 \in M$ not containing a critical point except possibly for their end points. 
By Lemma 1.8, $P(T(M))$ is the union of the regulated arcs $[P(z_1),P(z_2)]_K$. 
As this set is regulated connected and contains all of $P(M)$, 
by definition this set equals $[P(M)]_K$.
\endofproof
\medskip
{\bf 1.11 Remark.} 
If $X \subset K(P)$ is arbitrary, the same argument shows that 
$P(T(X)) \subset [P(X \cup \Omega(P))]_K$.
\bigskip
{\bf 1.12 Definition.} 
Let $T^*(M)$ be the family whose elements are the closures of 
components of $T(M)-\Omega(P)$. 
\bigskip
{\bf 1.13 Lemma.}
{\it P induces a continuous map from $T(M)$ to itself, 
where the restriction to every element (component) of $T^*(M)$ is injective.}
\medskip
{\bf Proof.} 
This follows from Lemmas 1.8 and 1.10.
\endofproof
\bigskip
{\bf 1.14 Lemma.}
{\it Let $\gamma(I) \subset K(P)$ be a regulated arc containing no critical value of $P$ except possibly for its end points. 
Then any lift of $\gamma(I)$ by $P$ is a regulated arc.}
\medskip
{\bf Proof.} 
As $\gamma|_{(0,1)}$ contains no critical value of $P$, 
it can be pulled back by $P$ in $d$ different ways, each being a regulated arc.
\endofproof
\bigskip
{\bf 1.15 Definition.} 
Given $z \in T(M)$ the {\it incidence number $\nu_{T(M)}(z)$ of $T(M)$ at $z$} is the number of components of $T(M)-\{z\}$. 
In other words, $\nu_{T(M)}(z)$ is the number of branches of $T(M)$ that are incident at $z$.
Note that this number might be different from the number of connected components of $K(P)-\{z\}$ 
(the incidence number at $z$ for $P$).
\smallskip
A point $z \in T(M)$ is called a {\it branching point of $T(M)$} if $\nu_{T(M)}(z)>2$ 
and an {\it end} if $\nu_{T(M)}(z)=1$. 
The {\it preferred set} of $T(M)$ is $V_{T(M)}= M \cup\break 
\{z \in T(M): \nu_{T(M)}(z)>2\}$. 
Note that $V_{T(M)}$ is finite. 
This because there are only a finite number of vertices in this tree.
\bigskip
{\bf 1.16 Proposition.}
{\it The set $V_{T(M)}$ is invariant. 
Furthermore, it generates the same tree as M; i.e, $T(M)=T(V_{T(M)})$.}
\medskip
{\bf Proof.} 
If $z$ is a branching point and $deg_z P=1$, 
then $P(z)$ is also a branching point with $\nu_{T(M)}(P(z)) \ge \nu_{T(M)}(z)$ because $P$ maps $T(M)$ into itself and $P$ is a local homeomorphism in a neighborhood of $z$.
\smallskip
We must prove that $[M]_K=[V_{T(M)}]_K$. 
As $M \subset V_{T(M)}$ then $[M]_K \subset [V_{T(M)}]_K$. 
Also by definition $V_{T(M)} \subset [M]_K$, 
so $[V_{T(M)}]_K \subset [[M]_K]_K=[M]_K$.
\endofproof
\bigskip
{\bf 1.17 Corollary.}
{\it Let $M, M'$ be finite invariant subsets containing $\Omega(P)$.  
If 
$V_{T(M)}=V_{T(M')}$ then $T(M)=T(M')$.}
\endofproof
\bigskip
{\bf 1.18 Proposition.}
{\it Let $v,v' \in J(P)$ be two periodic points. 
If for all $n \ge 0$, 
$P^{\circ n}(z)$ and $P^{\circ n}(z')$ belong to the same element (component) of $T^*(M)$, 
then $v=v'$.}
\medskip
{\bf Proof.} 
Suppose $P^{\circ n}(v),P^{\circ n}(v')$ belong to the same component of $T^*(M)$ for all $n \ge 0$. 
By Lemma 1.8 there is no precritical point in $[v,v']_K$. 
It follows easily that $[v,v']_K \subset J(P)$. 
Next, let $m$ be the least common multiple of the periods of $v$ and $v'$. 
Thus, $v,v'$ are fixed by $P^{\circ m}$. 
As there are only a finite number of such fixed points, we may assume that there are no other in this set $[v,v']_K$. 
Both endpoints of this regulated arc are repelling. 
Also by Lemma 1.8, $P^{\circ m}$ induces an homeomorphism of $[v,v']_K$ onto itself. 
It follows that there must be other fixed point in the interior of the arc $[v,v']_K$, 
in contradiction to what was assumed. 
\endofproof 
\bigskip
{\bf 1.19 Remark.} 
Note that the same is true if $v,v'$ are assumed only to be preperiodic. 
In this case, high enough iterates of both points must be periodic and therefore coincide. 
Lemma 1.13 will imply that $v,v'$ are identified as well. 
\bigskip
{\bf 1.20 Definition.}
We define the distance $d_{T(M)}(v,v')$ between points $v,v' \in V_{T(M)}$ as follows. 
Set $d_{T(M)}(v,v)=0$. 
Otherwise, take a regulated arc $[v,v']_K$ and define $d_{T(M)}(v,v')=\#(V_{T(M)} \cap [v,v']_K) - 1$ 
($\#$ denotes as usual cardinality). 
Thus, $d_{T(M)}$ measures the number of `edges' between $v$ and $v'$. 
In this language Proposition 1.18 can be read as follows. 
\bigskip
{\bf 1.21 Proposition: Expanding Property of the tree T(M).} 
{\it For all pairs
$v,v' \in V_{T(M)} \cap J(P)$ satisfying $d_{T(M)}(v,v')=1$, 
there is an $n \ge 1$ such that 
$d_{T(M)}(P^{\circ n}(v),P^{\circ n}(v'))>1$.} 
\medskip
{\bf Proof.} 
As $v,v'$ are eventually periodic, the result follows from Proposition 1.18.  
\endofproof
\insec
\centerline{\medfont 2. The Regulated Trees $T(P^{-n}M)$}
\insec
{In this section we study the inverse under $P$ of the tree $T(M)$.} 
\smallskip
{\bf 2.1 Proposition.}
{\it $P^{-1}T(M)=T(P^{-1}M)=T(P^{-1}V_{T(M)})$. 
In this case the vertices of the tree are given by 
$V_{T(P^{-1}M)}=P^{-1}V_{T(M)}$.}
\medskip
{\bf Proof.} 
As $P^{-1}M \subset P^{-1}V_{T(M)}$ we have $T(P^{-1}M) \subset T(P^{-1}V_{T(M)})$.
\smallskip
From Lemma 1.10, $PT(P^{-1}V_{T(M)})=[PP^{-1}V_{T(M)}]_K=[V_{T(M)}]_K=T(M)$. 
It follows that $T(P^{-1}V_{T(M)}) \subset P^{-1}T(M)$.
\smallskip
Now let $z \in P^{-1}T(M)-P^{-1}M$, then $P(z)$ belongs to a regulated arc $\gamma(I) \subset T(M)$, with only end points in $M$. 
By Lemma 1.14 any inverse of this regulated arc is also regulated with endpoints in $P^{-1}M$ and therefore belongs to $T(P^{-1}M)$; in this way 
$z \in T(P^{-1}M)$. 
If $z \in  P^{-1}M$ then by definition $z \in  T(P^{-1}M)$. 
This completes the proof of the chain of inequalities.
\smallskip
The second part follows from the first together with the definition of $V_{T(P^{-1}M)}$ and Proposition 1.16.
\endofproof
\bigskip
Proposition 2.1 extends easily.
\medskip
{\bf 2.2 Corollary.}
{\it $P^{-n}T(M)=T(P^{-n}M)=T(P^{-n}V_{T(M)})$. 
In this case the vertices of the tree are given by 
$V_{T(P^{-n}M)}=P^{-n}V_{T(M)}$.}
\endofproof
\bigskip
{\bf 2.3} 
As $T(M) \subset P^{-1}T(M)$ there are two incidence functions, 
$\nu_{0,M}=\nu_{T(M)}$ for $T(M)$ and $\nu_{-1,M}=\nu_{T(P^{-1}M)}$ for $P^{-1}T(M)$. 
It is immediate that $\nu_{0,M}(z) \le \nu_{-1,M}(z)$ at every point of $T(M)$. 
Furthermore, we have the following (here $deg_z P$ denotes the local degree of $P$ at $z$). 
\bigskip
{\bf 2.4 Lemma.}
{\it $\nu_{-1,M}(z)=\nu_{0,M}(P(z)) \hskip 0.05in deg_z P$, 
for every $z \in P^{-1}T(M)$.}
\medskip
{\bf Proof.} 
This follows from Lemma 1.10 and Proposition 2.1.  
\endofproof
\bigskip
These inequalities can be easily generalized for the incidence functions $\nu_{-n,M}$ of the trees $P^{-n}T(M)$. 
For example, $\nu_{-n,M}(z) \le \nu_{-n-1,M}(z)$ at every point of $P^{-n}T(M)$.
\bigskip
The next proposition is a weak attempt to reconstruct the tree $P^{-1}T(M)$ 
starting from $T(M)$. 
An improved version will be given in Chapter III 
(compare Proposition III.2.5). 
\bigskip
{\bf 2.5 Proposition.}
{\it Let $X$ be a component of $T^*(P^{-1}M)$. 
Denote by $C(X)=\Omega(P) \cap X$ the critical points in $X$. 
Then $P$ induces a homeomorphism between $X$ and the component $T_\alpha$ of $T(M)$ cut along $P(C(X))$ that contains $P(X)$.}
\medskip
{\bf Proof.} 
By Lemma 1.13 $P|X$ is injective. 
Also, $P(X)$ is relatively open in $T_\alpha$. 
As it is also compact it must be the whole component.
\endofproof
\eject
\insec 
\centerline{\medfont 3. Incidence.}
\insec 
{In this Section we take a closer look at terminal, incidence, branching and non 
branching points of the Postcritically Finite Polynomial $P$. 
A point $p \in J(P)$ is 
terminal if there is only one external ray landing at $p$. 
Otherwise $p$ is an incidence point. 
For incidence points we distinguish between branching (if there are more than two rays landing at $p$) 
and non branching (exactly two rays landing at $p$). 
We will show that for a postcritically finite polynomial $P$, 
every branching point must be periodic or preperiodic. 
Also we will prove that every periodic branching point is present as a preferred point (see $\S$1.15) in the minimal tree $T(M_0)$, 
and thus in any tree $T(M)$. 
} 
\bigskip
{\bf 3.1} 
Let $P$ be a Postcritically Finite Polynomial, and $z$ an arbitrary point in the Julia set $J(P)$. 
Every component of $J(P)-\{z\}$ is eventually mapped onto the whole Julia set, 
and therefore contains points whose orbit contains any specified point. 
We will use this fact in the following two propositions. 
\bigskip
{\bf 3.2 Proposition.}
{\it Let P be a Postcritically Finite Polynomial and\break 
$z \in J(P)$ a branching point. 
Then $z$ is preperiodic (or periodic).}
\medskip
{\bf Proof.}
Suppose $z$ does not eventually map to ${\cal O}(\Omega(P))$ 
(otherwise $z$ is already preperiodic). 
Fix $w \in \Omega(P)$ and pick in every component of $J(P)-\{z\}$ a point $p_i$ 
which eventually maps to $w$. 
The orbit ${\cal O}(\{p_1,\dots,p_k\})$ of this set $\{p_1,\dots,p_k\}$ is a finite set. 
In this way, the set $M'=M \cup {\cal O}(\{p_1,\dots,p_k\})$ is invariant and contains the critical points of $P$. 
As $z \in V_{T(M')}$, the result follows from Proposition 1.16.  
\endofproof
\bigskip
{\bf 3.3 Proposition.}
{\it Let P be a Postcritically Finite Polynomial and $z \in J(P)$ a 
periodic incidence point. 
Then $z \in T(M_0)$, and in this way $z \in T(M)$ for any finite 
invariant set $M \supset \Omega(P)$. 
Furthermore, $\nu_{0,M}(z)$ is independent of M and equals the 
number of components of $J(P)-\{z\}$. 
In particular there are exactly $\nu_{0,M_0}(z)$ external rays landing at z.}
\medskip
{\bf Proof.}
The number of rays landing at $z$ equals the number of components of 
$J(P)-\{z\}$. 
After this remark the proof is analogous to that of last proposition. 
Further details are left to the reader 
(compare also Lemma 1.10 and Remark 1.11). 
\endofproof
\bigskip
{\bf 3.4 Corollary.}
{\it Let $z \in J(P) \cap T(M)$ be such that $P^{\circ n}(z)$ is periodic. 
Then $\nu_{-n,M_0}(z)$ equals the number of components of $J(P)-\{z\}$.
In particular there are exactly $\nu_{-n,M}(z)$ external rays landing at z.}
\medskip
{\bf Proof.} 
This follows from Proposition 3.3 and Lemma 2.4. 
\endofproof
\bigskip
{\bf 3.5 Corollary.}
{\it T(M) contains a fixed point of P.}
\medskip
{\bf Proof.}
If $P$ has a fixed critical point, then such point is in $M$ and by definition in $T(M)$. 
Otherwise, as there are only $d-1$ fixed rays, but $d$ fixed points, 
one must be an incidence point. 
By Proposition 3.3, this fixed point is in $T(M)$.
\endofproof
\bigskip

%% file: chap.6.tex
\input psfig
\inchap
\centerline{\bigfont Chapter II} 
\centerline{\bigfont Abstract Hubbard Trees.}
\inchap
{In this Chapter we set our basic abstract framework. 
We carefully justify the importance of all the elements in the definition of abstract Hubbard Trees 
given in the introduction 
(compare Examples 2.11-13). 
In Section 1, we introduce some basic notation related to finite topological trees. 
In Section 2, we introduce dynamics in finite topological trees, 
and explain why further structure should be added in order to have a characterization of postcritically finite polynomials. 
In Section 3, 
the elements needed for this characterization are defined. 
In Section 4, we give a normalization in order to simplify notation, 
and we state our main result, 
namely necessary and sufficient conditions for the realization of Hubbard Trees. }
\insec 
\centerline{\medfont 1. Cyclic Trees.}
\insec 
{In this Section we only introduce some notation related to 
finite topological trees 
which would be used throughout the rest of this work. } 
\bigskip
{\bf 1.1 Definition.} 
By a {\it topological tree T} will be meant a finite connected acyclic 
$m$-dimensional simplicial complex ($m=0,1$). 
Given $p \in T$ we define the {\it incidence number $\nu_T(p)$ of $T$ at p} 
as the number of connected components of $T-\{p\}$. 
We say that $p \in T$ is a {\it branching point} if $\nu_T(p)>2$, 
and an {\it end} if $\nu_T(p)=1$. 
\medskip
A homeomorphism $\gamma:I=[0,1] \to T$ is called a {\it regulated path in T}. 
In general we will not distinguish between the map $\gamma$ and its image $\gamma(I)$. 
This because given two points $p,p' \in T$, 
any regulated path joining them will have the same image, 
which we denote by $[p,p']_T$. 
Given $X \subset T$ we denote by $[X]_T$ the smallest subtree of $T$ which contains $X$. 
Clearly this notation is compatible with that introduced before. 
\medskip
{\bf 1.2 Definition.} 
A {\it cyclic tree} is a triple $(T,V,\chi)$, where
\smallskip
(a) $T$ is the {\it underlying topological tree};
\tskip
(b) $V \subset T$ is finite set of {\it vertices} so that each component of $T-V$ is an open $1$-cell (an {\it edge});
\tskip
(c) For each $v \in V$, $\chi_v$ represents a {\it cyclic order} in the set 
$E_v=\{\ell_1,\dots,\ell_k\}$ of all edges with $v$ as a common endpoint.
\medskip
The presence of these $\chi_v$ 
naturally determines an isotopy class of embeddings of this tree $T$ into ${\bf C}$. \bigskip
{\bf 1.3 Pseudoaccesses.} 
If $\ell,\ell' \in E_v$ are consecutive in the cyclic order of $E_v$, 
we say that $(v,\ell,\ell')$ is a {\it pseudoaccess to $v$}. 
Take a pseudoaccess $(v,\ell,\ell')$ to $v$, 
and let the end points of the edge $\ell'$ be $v,v' \in V$. 
At $E_{v'}$ let $\ell''$ be the successor of $\ell'$ in the cyclic order. 
We say that $(v',\ell',\ell'')$ is the {\it successor} of $(v,\ell,\ell')$. 
\bigskip
{\bf 1.4 Lemma.}
{\it Let $(T,V,\chi)$ be a cyclic tree. 
The successor function in the set of pseudoaccess to the vertices in $V$ is a complete cyclic order.}
\smallskip
{\bf Proof.}
A trivial induction in the cardinality of $V$.
\endofproof
\bigskip
{\bf 1.5 Remark.} 
A Postcritically Finite Polynomial $P$ and a finite invariant set $M$ 
containing the critical set $\Omega(P)$ of $P$, 
naturally defines a cyclic tree $(T(M),V_{T(M)},\chi)$. 
Here $\chi_v$ represents the cyclic order of the components around a point $v \in V_{T(M)}$ taken counterclockwise. 
\bigskip
{\bf 1.6 Definition.} 
Let $(T,V,\chi)$ be a cyclic tree, and let $M \subset V$. 
We define the {\it restriction of $(T,V,\chi)$ to $M$}, 
as the cyclic tree
$([M]_T,V_M,\chi')$ 
where 
$V_M$ is the union of $M$ and the branching points in the topological tree $[M]_T$, 
and $\chi'_v$ is the natural restriction of the cyclic order $\chi_v$ of $E_v$ to 
the set $E'_v$ of all edges of $[M]_T$ incident at $v$. 
\insec 
\centerline{\medfont 2. Dynamical Abstract Trees.} 
\insec
{In this Section we give our first attempt to describe the dynamics of a Postcritically Finite polynomial by means of 
the dynamics in a finite topological tree. 
Unfortunately this simple characterization proves to be weak 
(compare Examples 2.11-13), and further structure has to be added. 
This will be done in Sections 3 and 4.} 
\bigskip
{\bf 2.1 Definition.} 
A {\it dynamical abstract tree} is a triple ${\bf T}=((T,V,\chi),\tau,\delta)$ where
\smallskip
(a) $(T,V,\chi)$ is {\it the underlying cyclic tree}, 
\tskip
(b) $\tau:V \to V$ is the {\it vertex dynamics},
\tskip
(c) $\delta:V \to {\bf Z}$ is a positive {\it local degree function}. 
\medskip
We require these elements to be related as follows, 
\smallskip
(i) For any edge $\ell$ with endpoints $v,v' \in V$ we must have $\tau(v) \ne \tau(v')$. 
\medskip
This condition allows us to extend $\tau$ to the underlying tree as follows. 
For any edge $\ell$ with endpoints $v,v' \in V$, 
map $\ell$ homeomorphically to the shortest path joining $\tau(v)$ and $\tau(v')$. 
Any extension `$\tau$' well defines a map $\tau_v:E_v \to E_{\tau(v)}$. 
We require, 
\smallskip
(ii) For any $v \in V$, there exists a cyclic ordered set ${\cal E}_v$ such that $E_v$ embeds in an order preserving way into ${\cal E}_v$. 
We require that $\tau_v$ can be extended to a degree $\delta(v)$ orientation preserving covering map between ${\cal E}_v$ and $E_{\tau(v)}$ 
(see appendix B). 
For the practical interpretation of this set ${\cal E}_v$ we refer to Remark 2.2 and Proposition III.2.5.
\medskip
We define the {\it degree} of ${\bf T}$ as 
$deg({\bf T})=1 + \sum_{v \in V}(\delta(v)-1)$. 
We require
\smallskip
(iii) $deg(T) > 1$. 
\bigskip
{\bf 2.2 Remark.} 
A Postcritically Finite Polynomial $P$ of degree $n>1$ and 
a finite invariant set $M \supset \Omega(P)$ 
naturally defines a dynamical abstract tree 
${\bf T}_{P,M}=((T(M),V_{T(M)},\chi),P,deg_z P)$ of degree $n$. 
Here ${\cal E}_v$ represents the components around $v$ in the tree $T^{-1}_M$ 
(see $\S$I.2). 
\bigskip
{\bf 2.3 Definitions.} 
Let ${\bf T}=((T,V,\chi),\tau,\delta)$ be a dynamical abstract tree. 
We extend $\delta$ to all the tree $T$ by letting $\delta(p)=1$ if $p \not\in V$. 
We define the {\it critical set of ${\bf T}$} as 
$\Omega({\bf T}) =\{p \in T: \delta(p)>1\}$. 
Condition (iii) above implies that $\Omega({\bf T})$ is always non empty. 
A point $p \in \Omega({\bf T})$ is a {\it critical point}; 
otherwise, it is {\it non critical}. 
\smallskip
The orbit of $S \subset V$ is the set ${\cal O}(S)=\cup^\infty_{k=0} \tau^{\circ k}(S)$. 
\bigskip
{\bf 2.4 Definition.} 
Let $\ell \in E_v$, 
we denote by ${\cal B}_{v,T}(\ell)$ the closure of the connected component of $T-\{v\}$ that contains $\ell$. 
This is just the {\it branch at $v$ determined by $\ell$ in the tree $T$}. 
\bigskip
{\bf 2.5 Definition.}
Let ${\bf T}=((T,V,\chi),\tau,\delta)$ be an abstract tree, 
and let $M \subset V$ be an invariant set of vertices containing the critical set 
($\tau(M) \cup \Omega({\bf T}) \subset M$). 
We define the {\it restriction ${\bf T}(M)$ of ${\bf T}$ determined by $M$}, 
as the abstract tree ${\bf T}(M)=((T(M),V_M,\chi),\tau',\delta')$, 
where 
$(T(M),V_M,\chi)$ is the restriction of the angled tree as defined in $\S$1.8, and  
$\tau'$, $\delta'$ are restrictions of the functions $\tau,\delta$ to the set $V_M$. 
\medskip
{\bf 2.6 Definition.}
Let ${\bf T,T'}$ be two abstract trees of degree $n=deg({\bf T})=deg({\bf T'}) > 1$. 
We say that ${\bf T'}$ is an {\it extension of ${\bf T}$} 
(in symbols ${\bf T} \preceq {\bf T'}$), 
if there is an embedding $\phi:T \to T'$ which satisfies the obvious conditions: 
\medskip
(i) $\phi(V) \subset V'$, 
\tskip
(ii) $\tau'(\phi(v))=\phi(\tau(v))$ and 
\tskip
(iii) $\delta(v)=\delta'(\phi(v))$ for all $v \in V$,
\tskip
(iv) $\phi$ induces a cyclic order preserving embedding of $E_v$ into $E_{\phi(v)}$. 
(At this point it is convenient to think of the elements of $E_v$ as `germs of edges'.) 
\medskip
Clearly $\preceq$ is an order relation.
\bigskip
{\bf 2.7} 
Let ${\bf T,T'}$ be two abstract trees of degree $n=deg({\bf T})=deg({\bf T'}) > 1$. 
We say that ${\bf T'}$ is {\it equivalent to ${\bf T}$} 
(in symbols ${\bf T} \approx {\bf T'}$),
if ${\bf T} \preceq {\bf T'}$ and ${\bf T'} \preceq {\bf T}$. 
This determines an equivalence relation between abstract trees. 
Furthermore, the order relation $\preceq$ extends to 
a partial order between equivalence classes of dynamical abstract trees of degree $n>1$. 
\bigskip
{\bf 2.8 Definition.} 
We say that an equivalence class ${\bf [T]}$ of dynamical abstract trees of degree $n>1$ is {\it minimal} if 
given ${\bf [T']} \preceq {\bf [T]}$ we necessarily have ${\bf [T']}={\bf [T]}$. 
\smallskip 
From the definition of extension tree we can deduce that 
if ${\bf [T']}$ is an extension of ${\bf [T]}$, 
then ${\bf [T]}$ is a restriction of ${\bf [T']}$ in the sense of Definition 2.5. 
Therefore we have the following. 
\bigskip
{\bf 2.9 Proposition.} 
{\it Every abstract tree ${\bf T}$ contains a unique minimal tree $min({\bf T})$.
Furthermore, this unique minimal tree is the tree generated by the orbit ${\cal O}(\Omega({\bf T}))$ of the critical set. }
\endofproof
\bigskip
{\bf 2.10}
The question now is if this description completely characterizes Postcritically Finite Polynomials. 
In other words, given a class ${\bf [T]}$ of dynamical abstract trees, is there a unique 
(up to affine conjugation) Postcritically Finite Polynomial $P$ and 
an invariant set $M \supset \Omega(P)$ such that ${\bf T}_{P,M} \in {\bf [T]}$? 
\smallskip
The answer is negative as the following examples show. 
\bigskip
{\bf 2.11 Non uniqueness.} 
Suppose a degree 3 polynomial has the following minimal tree ${\bf T}$ 
(where the double star stands for a double critical point, 
i.e, its local degree is 3). 
\bigskip
{\bigskip
\centerline{
\hbox to 4in{ 
\hbox to 2in{$\bullet \atop x_1$ {\hrulefill} $\star\star \atop x_0$} 
{\hrulefill}
$\bullet \atop x_2=x_3$}
}
\bigskip}

\bigskip
\centerline{\sl Figure 2.1. 
The vertex dynamics is given by $x_0 \mapsto x_1 \mapsto x_2 \mapsto x_3=x_2$.} 
\bigskip
If we want a centered monic polynomial with this minimal tree we suppose that $x_0=0$. 
We have then $P(z)=z^3+c$. 
(For polynomials of the form $P_c(z)=z^3+c$, the number $c^2$ is a complete invariant up to conjugacy. 
In other words $P_c$ is affine conjugate to $P_{c'}$ if and only if $c^2=c'^2$.)
If $P$ has this minimal tree, then the orbit of the critical point is as follows,
$0 \mapsto c \mapsto c^3+c \mapsto c^3+c$.
\smallskip
In this way, the relation $P_c^{\circ 2}(0)=P_c^{\circ 3}(0)$ determines the equation $c^3+c=(c^3+c)^3+c$. 
Thus $c$ must satisfy $c^5(c^4+3c^2+3)=0$. 
If we want $c^3+c \ne 0$ we must have $c \ne 0$, and we have two different possible values for $c^2={-3 \pm \sqrt{-3} \over 2}$. 
For both values of $c^2$ the respective polynomials $P_c$ have minimal tree $T(M_0)$ as shown in Figure 2.1. 
In fact, by Lemma I.1.13,\break 
$c$ and $c^3+c$ belong to different components of $T-\{0\}$. 
\smallskip
In this way, we have constructed two different non affine conjugate polynomials $P,P'$ which define the same class of minimal trees. 
Nevertheless, the trees ${\bf T}_{P,{\cal O}P^{-1}\Omega(P)}$, ${\bf T}_{P',{\cal O}P'^{-1}\Omega(P')}$ belong to different classes 
(see Figure 2.2 below). 
\bigskip

$$
\matrix{
& & &\ \bullet \ x_2^{-1}\hfill &&\hskip 0.2in 
&{{x_1 \hfill} \atop {\bullet {\hbox to 0.3in{\ \hrulefill} }}} \hskip -0.1in & {{x_0^{-1} \hfill}\atop {\bullet {\hbox to 0.3in{\ \hrulefill} }} } \hskip -0.1in
& {{x_0 \hfill} \atop {{\star\star} \ {\hbox to 0.3in{\hrulefill} }}} \hskip -0.1in
&{{x_0^{-1} \hfill} \atop {\bullet {\hbox to 0.3in{\ \hrulefill} }}} \hskip -0.2in
&  { x_2=x_3 \atop \bullet} \cr 
& & &\ \big\vert \hfill   && && &\ \big\vert \hfill       &\cr 
& & &\ \bullet \ x_0^{-1}\hfill & &&& &\ \bullet \ x_0^{-1}\hfill &\cr 
& & &\ \big\vert \hfill    &&& & &\ \big\vert \hfill      &\cr 
&{{\bullet {\hbox to 0.3in{\ \hrulefill} }} \atop {x_1 \hfill}} \hskip -0.1in
&{{\bullet {\hbox to 0.3in{\ \hrulefill} }} \atop {x_0^{-1} \hfill}} \hskip -0.1in
&{{{\star\star} \ {\hbox to 0.3in{\hrulefill} }} \atop {x_0 \hfill}} \hskip -0.1in 
&{{\bullet {\hbox to 0.3in{\ \hrulefill} }} \atop {x_0^{-1} \hfill}} \hskip -0.2in
& {\bullet \atop x_2=x_3} & & &\ \bullet \ x_2^{-1} &\cr
} 
$$
\bigskip
{\sl Figure 2.2. 
Here $x^{-1}_j$ maps to $x_j$. 
Even if the trees are isomorphic, they fail to have the same cyclic order around $x_0$.}
\bigskip
{\bf 2.12 Non existence} (compare Figure 2.3). 
The class of the tree below can not be obtained from a polynomial map. 
It must correspond to a degree two polynomial with three fixed points, 
which is impossible. 
{\bigskip
\centerline{
\hbox to 4in{ 
\hbox to 2in{$\bullet \atop y_0=y_1$ {\hrulefill} $* \atop x_0=x_1$} 
{\hrulefill}
$\bullet \atop z_0=z_1$}
}
\bigskip}
{\sl Figure 2.3.
All vertices are fixed. Here $\delta(x_0)=2$, and $\delta(y_0)=\delta(z_0)=1$.} 
\bigskip
Here is an alternative description of the obstruction for `realizing' this tree. 
If this tree is equivalent to a tree ${\bf T}_{P,M}$ of a degree two polynomial $P$, 
edges whose common vertex is the Fatou critical point 
must be realized near this vertex as internal rays in the uniformizing coordinate (see Section I.1.1). 
Let $\alpha$ be the difference of the two arguments in this coordinate. 
We have then $2\alpha=\alpha+1$ ($mod$ 1). 
But this implies that the two segments should be identified. 
Note that the minimal tree corresponding to this tree has only ${x_0}$ as vertex. 
Thus, this minimal tree can be realized as ${\bf T}_{z \mapsto z^2,\{0\}}$. 
\bigskip
{\bf 2.13 Non existence} (compare Figure 2.4). 
The class of the minimal tree below can not be obtained from a polynomial map. 
If there is a tree ${\bf T}_{P,M}$ in the class of such tree, 
it will not satisfy the expanding property 
(compare Propositions I.1.18 and I.1.21). 
\smallskip
\centerline{\psfig{figure=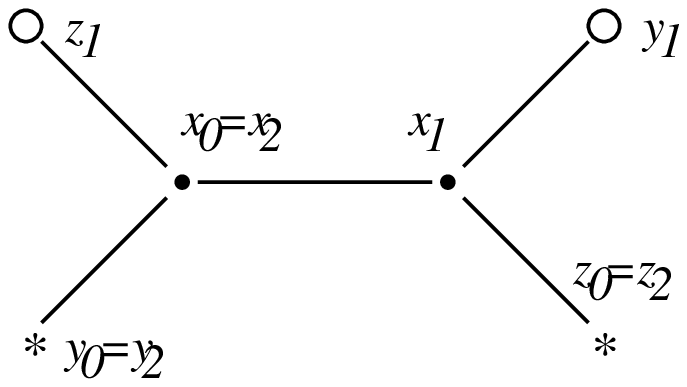,height=1.6in}} 
\centerline{\sl Figure 2.4. 
For any $k \ge 0$ there is no vertex between $\tau^{\circ k}(x_0)$ and $\tau^{\circ k}(x_1)$.} 
\bigskip
{\bf 2.14} 
All that can go wrong already happened in these three examples. 
Uniqueness failed because we had too little information. 
Here {\it too little information} means that we do not have enough information to recover in a unique way the tree $T^{-1}(M)$ 
(see Section I.2, compare also Propositions I.2.5 and III.2.5). 
Examples 2.12 and 2.13 failed because they do not satisfied necessary conditions. 
Namely, the trees must have well defined angles around Fatou critical points (see Section I.1) 
and should satisfy the expanding condition between Julia type vertices (see Proposition I.1.21). 
\insec 
\centerline{\medfont 3. Angled Trees.}
\insec 
{In this section we introduce the class of trees that will model our results. 
In this class we must be able to consider an analogue of the expanding condition, 
and also to define angles between edges near Fatou points.}  
\bigskip
{\bf 3.1 Definition.} 
An {\it angled tree} is a pair ${\bf A}=({\bf T},\angle)$, where
\smallskip
(a) ${\bf T}=((T,V,\chi),\tau,\delta)$ is a dynamical abstract tree, 
\tskip
(b) together with a function $\ell,\ell' \mapsto \angle(\ell,\ell')=\angle_v(\ell,\ell') \in {\bf Q/Z}$ 
which assigns a rational modulo 1 to each pair of edges $\ell,\ell'$ which meet at a common vertex. 
This angle $\angle(\ell,\ell')$ should be skew-symmetric, 
with $\angle_v(\ell,\ell')=0$ if and only if $\ell=\ell'$, 
and with 
$\angle_v(\ell,\ell'')=\angle_v(\ell,\ell')+\angle_v(\ell',\ell'')$ whenever three edges are incident at a vertex $v$. 
\medskip
The maps $\angle$, $\tau$ and $\delta$ must be related as follows. 
Again we extend $\tau$ to a map $\tau:T \to T$ 
which carries each edge homeomorphically onto the shortest path joining the images of its endpoints. 
Any extension well defines a map between `germs' $\tau_v:E_v \to E_{\tau(v)}$. 
We require then that 
$$
\angle_{\tau(v)}\big(\tau_v(\ell),\tau_v(\ell')\big)=\delta(v) \angle_{v}\big(\ell,\ell'\big), \eqno(1)
$$
whenever $\ell,\ell'\in E_v$ 
(in this case $\tau_v(\ell),\tau_v(\ell')$ contain edges incident at $\tau(v)$ where the angle between them is measured).
\medskip
Such an angle function determines a cyclic order in $E_v$ which we suppose to coincide with $\chi$. 
Note that in this case the angle function $\angle_v$ at $v$ can be extended to a bigger set ${\cal E}_v$ (see $\S$3.3 below). 
\medskip
The degree $deg({\bf A})$ of the angled tree ${\bf A}=({\bf T},\angle)$ is by definition the degree of the abstract tree ${\bf T}$. 
The critical set $\Omega({\bf A})$ of ${\bf A}$ is by definition $\Omega({\bf T})$. 
\smallskip
{\bf 3.2} 
A vertex $v \in V$ is called {\it periodic }if for some $m >0$ we have $\tau^{\circ m}(v)=v$. 
The orbit of a periodic critical vertex is a {\it critical cycle}. 
We say that a vertex $v$ is of {\it Fatou type (or a Fatou vertex)} if eventually maps to a critical cycle. 
Otherwise it is of {\it Julia type (or a Julia vertex)}. 
If $v_0 \mapsto v_1 \mapsto \dots \mapsto v_m=v_0$ is a critical cycle, 
we define the {\it degree of the cycle} as the product 
$\delta(v_0)\times \dots \times \delta(v_{n-1})$ 
of the degrees of the elements in said cycle. 
\bigskip
{\bf 3.3} 
The function $\tau$ induces a function $\tau_v$ between the set $E_v$ of edges incident at $v$ 
and the set $E_{\tau(v)}$ of edges incident at $\tau(v)$. 
Given a Fatou periodic vertex 
we can find embeddings $\phi_v \to {\bf R/Z}$ called {\it local coordinates of the set $E_v$} (see Appendix B) such that the diagram 
$$
\matrix{
&& E_v &{\tau_v \atop \longrightarrow} &E_{\tau(v)} &\cr 
&\phi_v &\Big\downarrow &&\Big\downarrow &\phi_{\tau(v)} \cr 
&& {\bf R/Z} &{m_v \atop \longrightarrow} &{\bf R/Z} &\cr} 
\eqno(2)
$$
commutes. 
Here $m_v$ is multiplication by $\delta(v)$ (modulo 1). 
Note that the number of possible embeddings for each critical cycle is the degree of the cycle minus one. 
\smallskip
At other Fatou vertices $v$ 
we can still make diagram (2) hold by pulling back the local coordinate at $\tau(v)$ 
and using relation (1). 
\smallskip
At periodic Julia vertices relation (1) easily implies that $\tau_v$ is a bijection. 
We pick an element $\ell \in E_v$ to which we assign the 0 coordinate 
($\phi_v(\ell)=0$). 
If $E_{\tau(v)}$ has not been assigned a local coordinate, 
we assign to each edge $\tau_v(\ell) \in E_{\tau(v)}$ the argument $\phi_v(\ell)$. 
In general we can not make diagram (2) commute for all vertices. 
(It might fail at the starting vertex $v$). 
In this last case the induced function $m_v$ in ${\bf R/Z}$ becomes translation by some constant.   
\smallskip
At non periodic Julia vertices a local coordinate $\phi_v$ can be pulled back from $E_{\tau(v)}$ in $\delta(v)$ different ways such that diagram (2) commutes.  
\bigskip
{\bf 3.4 Definition.}
Let ${\bf A}=({\bf T}=((T,V,\chi),\tau,\delta),\angle)$ be an angled tree. 
For a finite set of invariant vertices $M \supset \Omega({\bf A})$, we denote by 
${\bf A}(M)=({\bf T}(M),\angle_M)$ the angled tree generated by $M$, 
i.e, take ${\bf T}(M)$ the dynamical abstract tree determined by $M$ (see section 2.5), 
and let $\angle_M$ be the restriction of $\angle$ to the vertices of ${\bf T}(M)$. 
Of course, ${\bf A}={\bf A}(V)$. 
\bigskip
{\bf 3.5 Lemma.}
{\it For any extension `$\tau$' and invariant set of vertices $M \supset \Omega({\bf A})$ we have $\tau(T(M))=[\tau(M)]_T$.}
\smallskip
{\bf Proof.}
A copy of Lemma I.1.10 with the appropriate change of notation (see also Lemma 3.7). 
\endofproof
\bigskip
{\bf 3.6 Definition.}
Let ${\bf A}$ be an angled tree of degree $n$, and 
let $\Omega({\bf A})=\{v_1,\dots,v_l\}$ be the critical set. 
For a fixed family of local coordinates $\{\phi\}_{v \in V}$, 
we construct a partition $T^*=T^*(\{\phi\})$ of $T$ consisting (counting possible repetitions) of exactly $n$ subtrees of $T$. 
This partition will have the property that every point $p \in T$ will belong to exactly $\delta(p)$ elements of $T^*$. 
Note that this will be possible only if we somehow `unglue' the tree around every critical point. 
This is formally done as follows 
(compare the example within the proof of Proposition III.2.5). 
\smallskip
Let $T_0$ be $\{T\}$. 
We will inductively define partitions $T_i$ ($i \le l$) of $T$ with the following properties 
\smallskip
(a) For $j \le i$, $v_j$ belongs to exactly $\delta(v_j)$ elements of $T_i$, 
\tskip
(b) For $j >i$, $v_j$ belongs to exactly one element of the family $T_i$,
\tskip
(c) $T_i$ is constructed from $T_{i-1}$ 
by replacing the unique element $T(\alpha)$ of $T_{i-1}$ to which $v_i$ belongs 
by $\delta(v_i)$ subtrees of $T(\alpha)$. 
\smallskip
We proceed as follows. Let $T(\alpha)$ be the only element of $T_{i-1}$ to which $v_i$ belongs. 
We partition $T(\alpha)$ into $\delta(v_i)$ pieces as follows. 
First divide the set $E_i=E_{v_i}$ in $\delta(v_i)$ subsets using the local coordinate. 
For this we define for $k=0,\dots, \delta(v_i)-1$,  
$$
E^k=E^k_i=\{\ell \in E_i: \phi_{v_i}(\ell) \in [{k \over \delta(v_i)},{k+1 \over \delta(v_i)})\}
$$
Now, we take the union of all branches in a set $E^k$, i.e, define
$$
T^k(\alpha)= T(\alpha) \cap \big( v_i \cup \bigcup_{\ell \in E^k} {\cal B}_{v_i,T}(\ell)\big). 
$$
Define now $T_i$ by removing $T(\alpha)$ of the family and including all such $T^k(\alpha)$. 
By definition $T^*$ is the last partition $T_l$. 
\bigskip
{\bf 3.7 Lemma.}
{\it Let ${\bf A}$ be an angled tree. 
Then the vertex dynamics $\tau$ induces a continuous map of T into itself, 
where the restriction to every element (component) of $T^*$ is injective.}
\smallskip
{\bf Proof.} 
(Compare Lemma I.1.8.)
Let $T_\alpha$ be an element of $T^*$. 
Suppose there are different $p_1,p_2 \in T_\alpha$ so that $\tau(p_1)=\tau(p_2)$. 
Take a path $\gamma:I \to [p_1,p_2]_T \subset T_\alpha$ joining $p_1,p_2$. 
As $\tau|_{T_\alpha}$ is locally one to one, 
the set $\Delta=\{(t_1,t_2): t_1<t_2$ and $\tau(\gamma(t_1))=\tau(\gamma(t_2))\}$ is compact. 
As we have assumed that this set is not empty we can take $(t_1,t_2) \in \Delta$ with $t_2-t_1$ minimal. 
Let $t \in (t_1,t_2)$, then $\tau(\gamma([t_1,t]))$ and $\tau(\gamma([t,t_2]))$, are regulated arcs with the same end points. 
Therefore they are equal and thus $t_2-t_1$ is not minimal.
\endofproof
\bigskip
{\bf 3.8 Remark.}
As $T^*$ consists of $n=deg({\bf A})$ elements (counting possible repetitions), 
it follows from the last lemma that for any $p \in \tau(T)$ and any possible extension `$\tau$'
$$ \sum_{\{q \in T: \tau(q)=p\}} \delta(q) \le n.$$
\bigskip
{\bf 3.9 Lemma.}
{\it Let $v$ be a periodic Fatou vertex, and $\ell_1,\ell_2 \in E_v$ be different edges. 
There is an $n \ge 0$ so that $\tau_v^{\circ n}(\ell_1)$ and $\tau_v^{\circ n}(\ell_2)$ belong to different components of $T^*$.}
\smallskip
{\bf Proof.} 
Let $d>1$ be the degree of the cycle $v_0=v \mapsto v_1 \mapsto \dots \mapsto v_m=v_0$. 
We write $\phi_v(\ell_1)$ and $\phi_v(\ell_2)$ in base $d$ expansion. 
If for all $n$, $\tau_v^{\circ n}(\ell_1)$ and $\tau_v^{\circ n}(\ell_2)$ belong to the same component of $T^*$, 
by construction for all $k>0$ the integer parts of 
$m_{\delta(v_k)}\phi_{v_k}(\tau^{\circ k}(\ell_1))$ and $m_{\delta(v_k)}\phi_{v_k}(\tau^{\circ k}(\ell_2))$ are equal. 
But this implies that $\ell_1=\ell_2$. 
\endofproof
\bigskip
{\bf 3.10 Definition.} (Compare $\S$I.1.20.) 
We define the distance $d_T(v,v')$ between vertices as follows. 
Set $d_T(v,v)=0$. 
Otherwise let 
$d_T(v,v')$ be the number of edges between $v$ and $v'$. 
\smallskip
We say that the angled tree ${\bf A}=({\bf T},\angle)$ is {\it expanding} if the following property is satisfied 
(see also Propositions I.1.18 and I.1.21). 
\medskip
{\it
For any edge $\ell$ whose end points $v,v'$ are Julia vertices there is an $m \ge 1$ such that $d_T(\tau^{\circ m}(v),\tau^{\circ m}(v')>1$.}
\medskip
Equivalently, ${\bf A}$ is not expanding if and only if 
{\it there exists periodic Julia vertices $v,v'$ 
such that $d_T(\tau^{\circ m}(v),\tau^{\circ m}(v'))=1$ for all $m \ge 0$.}
\bigskip
{\bf 3.11 Lemma.}
{\it An angled tree ${\bf A}$ is expanding 
if and only if 
for any two periodic Julia vertices $v,v'$ there is an $m \ge 0$ such that 
$\tau^{\circ m}(v)$ and $\tau^{\circ m}(v')$ belong to different components of $T^*$.}
\medskip
{\bf Proof.}
Suppose ${\bf A}$ is not expanding. 
By definition there are periodic Julia vertices $v,v'$ with $d_T(\tau^{\circ m}(v),\tau^{\circ m}(v'))=1$ for all $m \ge 0$. 
As there are no critical points in the orbit of periodic Julia vertices, 
by construction $\tau^{\circ m}(v)$ and $\tau^{\circ m}(v')$ will be in the same element of $T^*$ for any possible choice of the family $\{\phi_v\}$. 
\smallskip
Let now ${\bf A}$ be expanding. 
Suppose there are different Julia vertices $v,v'$ such that 
$\tau^{\circ m}(v)$,$\tau^{\circ m}(v')$ belong to the same component of $T^*$ for all $m \ge 0$. 
Among such pairs we can take $v,v'$ periodic and with the property that $d_T(v,v')$ is minimal. 
By assumption the regulated path $[\tau^{\circ m}(v)\tau^{\circ m}(v')]_T$ is completely contained within a component of $T^*$ for all $m \ge 0$. 
It follows from Lemma 3.7 that all $[\tau^{\circ m}(v)\tau^{\circ m}(v')]_T$ are homeomorphic. 
We take $v'' \in [v,v']_T \cap V$ such that $d_T(v,v'')=1$. 
As ${\bf A}$ is expanding it follows that $v''$ is a periodic Fatou vertex. 
In this way $E_{v''} \cap [v,v']_T=\{\ell_1,\ell_2\}$ with $\ell_1 \ne \ell_2$. 
We get a contradiction in applying Lemma 3.9. 
\endofproof
\bigskip
{\bf 3.12 Corollary.} 
{\it Let ${\bf A}$ be an expanding angled tree. 
The induced angled tree ${\bf A}(M)$ is expanding for every invariant set of vertices $M \supset \Omega({\bf A})$.}
\endofproof
\bigskip
{\bf 3.13 Lemma.} 
{\it Let ${\bf A}$ be an expanding angled tree. 
Given a periodic Julia vertex $v$, every component of $T-\{v\}$ contains a vertex 
which belongs to ${\cal O}(\Omega({\bf A}))$.}
\smallskip
{\bf Proof.}
Suppose that ${\cal B}_{v,T}(\ell)$ does not contain a vertex in ${\cal O}(\Omega({\bf A}))$ different from $v$ for some $\ell \in E_v$. 
The relation $\tau_v(\ell) \in E_{\tau(v)}$ determines a cyclic sequence of edges 
$\ell=\ell_0 \in E_v, \ell_1 \in E_{\tau(v)},\dots,\ell_n=\ell_0 \in E_{\tau^{\circ m}(v)}=E_v$. 
If for some $k < m$ the branch ${\cal B}_{\tau^{\circ k}(v),T}(\ell_k)$ contains a critical point, 
we may assume that $k$ is as big as possible and derive a contradiction by using Lemma 3.7. 
We assume though that ${\cal B}_{\tau^{\circ k}(v),T}(\ell_k)$ does not contains a critical point for all $k$. 
This implies using again Lemma 3.7 that all ${\cal B}_{\tau^{\circ k}(v),T}(\ell_k)$ are homeomorphic with only periodic Julia vertices. 
Thus, ${\bf A}$ is not expanding.
\endofproof
\insec 
\centerline{\medfont 4. Abstract Hubbard Trees.}
\insec 
{The angles at Julia vertices are rather artificial, so we normalize them as follows. 
If $m$ edges $\ell_1,\dots,\ell_m$, meet at a periodic Julia vertex $v$, 
then we assume that the angles $\angle_v(\ell_l,\ell_k)$ are all multiples of $1/m$ 
(it follows that the angles at periodic Julia vertices convey no information beyond the cyclic order of these $m$ incident edges). 
Fortunately, 
this number is preserved under restrictions which contain the orbit of the critical set. 
This will allow us to give a coherent description.} 
\bigskip
{\bf 4.1 Definition.}
By an {\it abstract Hubbard Tree} we mean an expanding angled tree ${\bf H}=({\bf T},\angle)$ such that 
the angles at any periodic Julia vertex where $m$ edges meet are multiples of $1/m$. 
\bigskip
{\bf 4.2} 
Let ${\bf H,H'}$ be two abstract Hubbard Trees of degree $n=deg({\bf H})=deg({\bf H'})>1$. 
We say that ${\bf H'}$ is an {\it extension} of ${\bf H}$
(in symbols ${\bf H} \preceq {\bf H'}$), 
if there is an embedding $\phi:T \to T'$ which satisfies the obvious conditions: 
\medskip
(i) $\phi(V) \subset V'$, 
\tskip
(ii) $\tau'(\phi(v))=\phi(\tau(v))$ and 
\tskip
(iii) $\delta(v)=\delta'(\phi(v))$ for all $v \in V$,
\tskip
(iv) $\angle_v(\ell,\ell')=\angle'_{\phi(v)}(\phi(\ell),\phi(\ell'))$ for all $\ell,\ell' \in E_v$. 
\medskip
Clearly $\preceq$ is an order relation.
\bigskip
{\bf 4.3} 
Let ${\bf H,H'}$ be two abstract Hubbard Trees of degree $n=deg({\bf H})=deg({\bf H'}) > 1$. 
We say that ${\bf H'}$ is {\it equivalent to ${\bf H}$} 
(in symbols ${\bf H} \cong {\bf H'}$),
if ${\bf H} \preceq {\bf H'}$ and ${\bf H'} \preceq {\bf H}$. 
\medskip
This determines an equivalence relation between abstract Hubbard\break 
Trees. 
Furthermore, the order relation $\preceq$ well defines a partial order between 
equivalence classes of abstract Hubbard Trees of degree $n>1$. 
\bigskip
{\bf 4.4 Lemma.}
{\it Let ${\bf H}$ be an abstract Hubbard Tree, and $M \supset \Omega({\bf H})$ a finite invariant set of vertices. 
Then ${\bf H}(M)$ is an abstract Hubbard Tree and ${\bf H}(M) \preceq {\bf H}$.}
\smallskip
{\bf Proof.} 
This follows from Corollary 3.12 and Lemma 3.13.
\endofproof
\bigskip
{\bf 4.5 Proposition.} 
{\it Every abstract Hubbard Tree ${\bf H}$ contains a unique minimal tree $min({\bf [H]})$.
Furthermore, this unique minimal tree is the tree generated by the orbit ${\cal O}(\Omega({\bf H}))$ of the critical set.}
\smallskip
{\bf Proof.} 
This follows from Proposition 2.8 and Lemma 4.4.
\endofproof
\bigskip
{\bf 4.6 Remark.} 
A Postcritically Finite Polynomial $P$ and a finite invariant set $M \supset \Omega(P)$ naturally defines an abstract Hubbard Tree ${\bf H}_{P,M}=({\bf T}_{P,M},\angle)$. 
To define the angle function we note the following. 
At Fatou periodic vertices the edges of the tree are by definition segments of constant argument in the B\"ottcher coordinate (see I.1.1),  we define the angle between two such edges as the difference of their coordinates. 
For other Fatou points the coordinate can be defined such that the diagram (2) commutes, and we proceed as above. 
For a Julia set point $v$, $J(P)-\{v\}$ consists of a finite number (say $m$) of components. 
We define the `angle' between these components to be a multiple of $1/m$. 
As edges in the tree correspond locally to some of these components we have an angle function between them. 
(This procedure is well defined and compatible with the definition above, see Proposition I.3.3). 
It is easy to see that the minimal tree $T_{M_0}$ (see I.1.9) corresponds to the minimal tree ${\bf H}_{P,M_0}=min({\bf H}_{P,M})$ of any bigger invariant set $M$. 
\medskip
The main result of this work is the following. 
\bigskip
{\bf 4.7 Theorem.} 
{\it Let ${\bf H}$ be an abstract Hubbard Tree. 
Then there is a unique (up to affine conjugation) Postcritically Finite polynomial P, and an invariant set $M \supset \Omega(P)$ such that ${\bf H}_{P,M} \in {\bf [H]}$.} 
\bigskip
{\bf 4.8 Theorem.} 
{\it Equivalence classes of minimal abstract Hubbard Trees of degree $n >1$ are in one to one correspondence with affine conjugate Post-critically finite polynomials.}
\medskip
We prove Theorem 4.7 in the next chapter. 
Theorem 4.8 is an easy consequence of this result and Proposition 4.5.
\bigskip 

%% file: chap7.tex
\inchap
\centerline{\bigfont Chapter III}
\centerline{\bigfont Realizing Abstract Hubbard Trees.}
\inchap
{In this chapter we give the proof of the realization Theorem for Abstract Hubbard Trees (Theorem II.4.7). 
Our proof depends in the theory of Critical Portraits developed in the first part of this work. 
In Section 1 we define the class of extensions which do not add any essential information to the tree. 
We will prove later that every extension belongs to this class 
(compare Corollary 4.6). 
Section 2 gives the abstract analogue of $\S$I.2, where we show that a Hubbard 
Tree contains all the information required to reconstruct its `inverse'. 
Section 3 gives the abstract analogue of $\S$I.3. 
In Section 4 we relate the `accesses to Julia points' with the argument of a 
possible `external ray' (compare Theorem B in the introduction). 
As a consequence of this, 
we prove that every extension of a Hubbard Tree 
is canonical in the sense described in Section 1. 
In Section 5 we associate a Formal Critical Portrait to our Tree. 
This Critical Portrait is also admissible as shown in Section 6. 
Finally we prove that the Hubbard Tree associated with this critical portrait is 
equivalent to the starting one, thus establishing the result. 
From now on, we omit the trivial case in which $T$ is a single critical vertex.} 
\insec 
\centerline{\medfont 1. Canonical Extensions.}
\insec
{In this Section we define what we call `canonical extensions'. 
We will prove in Section 4 that every extension which itself is a Hubbard Tree, 
is canonical in the sense described here. 
This fact will allow us later to associate in a natural way a critical portrait to every Hubbard Tree.} 
\bigskip
{\bf 1.1 Definition.}
Let ${\bf H}_0 \preceq {\bf H}_1$ be abstract Hubbard Trees. 
We say that ${\bf H_1}$ is a {\it canonical extension} of ${\bf H}_0$ if 
for every extension ${\bf H} \succeq {\bf H}_0$, 
there is a common extension of ${\bf H}$ and ${\bf H}_1$. 
Canonical extensions always exist. 
By definition every Hubbard Tree is a canonical extension of itself. 
Our final goal in this direction will be to prove that every extension is canonical
(compare Corollary 4.6).
\bigskip
{\bf 1.2 Proposition.} 
{\it Let ${\bf H}$ be an abstract Hubbard Tree and $\omega$ a periodic Fatou vertex. 
There is a canonical extension ${\bf H'}$ of ${\bf H}$ such that 
\smallskip
(a) $E_v=E'_v$ at all vertices of the original Hubbard Tree ${\bf H}$. 
\tskip
(b) For every periodic $\ell \in E'_\omega$ with end points $\omega,v$ in ${\bf H}'$, the vertex $v$ is of Julia type and 
$d_{{\bf H}'}(\tau^{\circ m}(\omega),\tau^{\circ m}(v))=1$ for all $m \ge 0$.
\smallskip
In fact, the underlying topological trees can be chosen to be the same, with  
only new Julia vertices to be added.} 
\medskip 
{\bf Proof.}
Suppose the edge $\ell$ has end points $\omega,v$, and its germ is of period $k$. 
In other words suppose  
the induced maps $\tau_v$ determine a periodic sequence of edges 
$\ell_0=\ell \in  E_\omega,  
\ell_1 \in E_{\tau(\omega)}, \dots, \ell_k=\ell_0 \in E_{\tau^{\circ k}(\omega)}=E_{\omega}$. 
We distinguish two cases. 
\smallskip
Suppose $d_{\bf H}(\tau^{\circ m}(\omega),\tau^{\circ m}(v))=1$ for all $m$. 
If $v$ is of Julia type, condition (b) is already satisfied. 
If $v$ is of Fatou type then by Lemma II.3.7 all $\ell_k=[\tau^{\circ m}(\omega),\tau^{\circ m}(v)]_T$ are homeomorphic. 
In this case we insert a vertex $v_m$ in each $\ell_m$ (if $\ell_m=\ell_l$ then $v_m=v_l$) 
and define $\tau(v_m)=v_{m+1}$. 
Then clearly $v_1$ is periodic of period $k$ or $k/2$. 
The angles at $v_k$ are $1/2$ because two edges will meet now. 
Note that in this case this is the only possible extension that involves the segments $[\tau^{\circ m}(\omega),\tau^{\circ m}(v)]_T$ 
and gives an expanding tree.  
\smallskip
Otherwise, suppose $d_{\bf H}(\tau^{\circ m}(\omega),\tau^{\circ m}(v))>1$ for some $m \ge 1$.
In this case we insert a vertex $v_m$ in each $\ell_m$ as close as possible to 
$\tau^{\circ m}(\omega)$
(note here that if $\ell_m=\ell_l$ then we must have $v_m \ne v_l$) 
and define $\tau(v_m)=v_{m+1}$. 
Clearly $v_1$ is periodic of period $k$. 
The angles at $v_k$ are $1/2$ because two edges meet now. 
\medskip
The only obstruction to this construction is if condition (b) is already satisfied. 
Therefore the extension is canonical. 
\endofproof
\bigskip
{\bf 1.3 Corollary.}
{\it Every abstract Hubbard Tree has a canonical extension with at least one Julia vertex.}
\endofproof
\insec 
\centerline{\medfont 2. Inverse Hubbard Trees.}
\insec 
{We now describe an important type of canonical extension. 
In the case of the  Hubbard Tree ${\bf H}_{P,M}$ generated by a polynomial $P$ and an invariant set $M$, 
the interpretation is simple. 
We will reconstruct the equivalence class of the abstract Hubbard Tree generated by $P^{-1}M$ starting from ${\bf H}_{P,M}$. 
Thus, this section is the abstract analogue of Section I.2.} 
\bigskip
{\bf 2.1 Definition.} 
An abstract Hubbard Tree ${\bf H}$ of degree $n>1$ is\break 
{\it homogeneous} if 
\tskip
(a) $\forall v \in \tau(V),  n=\sum_{\{v' \in V: v=\tau(v')\}} \delta(v')$, and  
\tskip
(b) $\Omega({\bf H}) \subset \tau(V)$. 
\smallskip
In other words, 
every vertex with at least one inverse must have a maximal number counting multiplicity 
(compare Remark II.3.8). 
Furthermore, all critical vertices must have a preimage. 
The terminology is justified by the fact (proved below) that 
the underlying topological tree can be `chopped' into $n$ pieces; 
each piece being homeomorphic as a graph to the abstract tree generated by restriction to $\tau(V)$. 
More formally, 
$\tau$ establishes a homeomorphism between each of the $n$ elements of $T^*$ (compare Section II.3.6) 
and the abstract Hubbard Tree ${\bf H}(\tau(V))$. 
\bigskip
{\bf 2.2 Lemma.} 
{\it For any election of local coordinate system $\{\phi_v\}_{v \in V}$, 
each $T_\alpha \in T^*(\{\phi_v\})$ is homeomorphic to ${\bf H}(\tau(V))$.} 
\medskip
{\bf Proof.}
By Lemma II.3.5 we have $\tau(T)=[\tau(V)]_T$. 
Also, every $v \in  [\tau(V)]_T$ has at most one inverse in 
$T_\alpha$ by Lemma II.3.7. 
It follows easily from condition a) 
that every $v \in  [\tau(V)]_T$ should have a unique inverse in $T_\alpha$. 
The result follows. 
\endofproof
\bigskip
{\bf 2.3 Corollary.}
{\it Let ${\bf H}$ be an abstract tree of degree $n>1$, 
such that $\Omega({\bf H}) \subset \tau(V)$. 
Then ${\bf H}$ is homogeneous if and only if 
$\#(V)-1=\break 
n(\#(\tau(V))-1)$.}
\medskip
{\bf Proof.} 
This follows from Lemma 2.2 and Remark II.3.8. 
\endofproof
\bigskip
{\bf 2.4 Definition.}
Let ${\bf H'} \preceq {\bf H}$ be abstract Hubbard Trees with ${\bf H}$ homogeneous. 
We say that ${\bf H'}$ is the image of ${\bf H}$ if 
the embedding which defines the order ${\bf H'} \preceq {\bf H}$ is such that also 
${\bf H}(\tau(V)) \cong {\bf H'}$. 
\smallskip
This definition clearly extends to equivalence classes of abstract Hubbard Trees. 
\bigskip
{\bf 2.5 Proposition.}
{\it Every equivalence class ${\bf [H]}$ of abstract Hubbard Trees is the image of a unique class of homogeneous abstract Hubbard Trees.}
\medskip
{\bf Proof.} 
The proof of existence is constructive using only necessary conditions, 
uniqueness follows. 
Let $\{\phi_v\}_{v \in V}$ be a family of local coordinates for $V \subset T$. 
We will work with the family $T^*=T^*(\{\phi_v\})$ (compare $\S$II.3.6). 
We construct a new simplicial complex by 
gluing a different copy of $T$ to each component $T_\alpha \in T^*$ following $\tau$
(compare Lemma 2.2). 
In other words we consider $n$ disjoint copies $H^\alpha$ of $T$ 
($\alpha=1 \dots n$), 
with a suitable identification at ``critical points" described below. 
By Lemma II.3.7, the dynamics $\tau$ restricted to each subtree $T_\alpha$ 
of the family $T^*=T^*(\{\phi_v\})$ is one to one. 
We denote this restriction by $i_\alpha$ ($``i"$ stands for identification). 
Thus we have a family of maps $i_\alpha:T_\alpha \to H^\alpha$. 
We establish an equivalence relation $\sim$ between points in 
the disjoint union $\coprod H^\alpha$ as follows. 
Whenever $\omega \in T_\alpha \cap T_\beta$ (and this can only happen if $\omega$ is critical), 
we write 
$i_\alpha(\omega) \sim i_\beta(\omega)$. 
Thus the new underlying topological tree is 
$X=\coprod H^\alpha/ \sim$. 
There is a `natural inclusion $T \subset X$' induced by the maps $i_\alpha$. 
The new set of vertices is the disjoint union of vertices of $H_\alpha$ modulo 
$\sim$. 
\medskip
In order to avoid confusion in the above notation, 
we will interrupt the proof in order to exemplify our construction. 
$$
{\bullet \atop {x_1}} \hbox to 0.7in{\hrulefill} {** \atop {x_0}}   
\hbox to 0.7in{\hrulefill}{* \atop {x_2=x_3}} 
$$
The abstract tree in the figure above can be chopped into $4$ pieces according to the construction in $\S$II.3.6. 
We think of these pieces as mapping onto different copies $H_\alpha$ of $T$ 
(this is emphasized below by the superscripts in the right). 
 
$$
\matrix{
&T_1: &{\bullet \atop {x_1}} \hbox to 0.7in{\hrulefill} {\bullet \atop {x_0}} 
&{i_1 \atop \longrightarrow} &
{\bullet \atop {x^1_1}} \hbox to 0.7in{\hrulefill} {\bullet \atop {x^1_0}}   
\hbox to 0.7in{\hrulefill}{\bullet \atop {x^1_2}} \cr\cr\cr 
&T_2: &{\bullet \atop {x_0}} 
&{i_2 \atop \longrightarrow} &
{\bullet \atop {x^2_1}} \hbox to 0.7in{\hrulefill} {\bullet \atop {x^2_0}}   
\hbox to 0.7in{\hrulefill}{\bullet \atop {x^2_2}} \cr \cr\cr
&T_3: &{\bullet \atop {x_0}}   
\hbox to 0.7in{\hrulefill}{\bullet \atop {x_2=x_3}}
&{i_3 \atop \longrightarrow} &
{\bullet \atop {x^3_1}} \hbox to 0.7in{\hrulefill} {\bullet \atop {x^3_0}}   
\hbox to 0.7in{\hrulefill}{\bullet \atop {x^3_2}} \cr\cr\cr 
&T_4: &{\bullet \atop {x_2}} 
&{i_4 \atop \longrightarrow} &
{\bullet \atop {x^4_1}} \hbox to 0.7in{\hrulefill} {\bullet \atop {x^4_0}}   
\hbox to 0.7in{\hrulefill}{\bullet \atop {x^4_2}}  
}
$$ 
\bigskip
In this way the new tree is given by identifying $x^1_1=i_1(x_0),x^2_1=i_2(x_0)$ and $x^3_1=i_3(x_0)$ (because $x_0 \in T_1 \cap T_2 \cap T_3$) 
and by identifying 
$x^3_2=i_3(x_2)$ with $x^4_2=i_4(x_2)$ (because $x_2 \in T_3 \cap T_4$). 
Note that the original tree is canonically embedded in this new one by using $i_\alpha$. 
\bigskip
{\bf Proof of 2.5 (Continue).}
We continue the proof by defining the dynamics and angle functions. 
What we have done so far is simply to replace each piece $T_\alpha$ by the copy $H^\alpha$. 
In this way, 
if we think of the $H^\alpha$ as the corresponding pieces for the new tree, 
the final structure is induced by the old one by gluing the $H^\alpha$ following that same pattern of the $T_\alpha$. 
\medskip
The vertex dynamics $\bar{\tau}$ maps each new vertex to the actual point in $V$ from which it was constructed. 
More formally, 
take $v \in H^\alpha$ a vertex of $X$; 
as $H^\alpha$ is also partitioned by the family ${\bf T^*}$, it follows that 
$v \in T_\beta$ 
for some $\beta$. 
We define 
$\bar{\tau}(v)=i_\beta(v) \in H_\beta \subset X$. 
(Clearly this is well defined and two consecutive vertices have different image). 
The degree is one at each vertex not present in the original tree. 
In other words, 
if $v \in T_\alpha$ (that is if $v$ belongs to the original tree $T$), 
we define the degree at $i_\alpha(v)$ 
(which is the point in $X$ to which $v$ is identified) 
as $\bar{\delta}(i_\alpha(v))=\delta(v)$.   
If $v \in X$ is not of the form $i_\beta(\omega)$ for some $\omega$, 
we set $\bar{\delta}(v)=1$. 
\smallskip
The angle function at non critical points is pulled back from the identification: 
if $\bar{\delta}(w)=1$, we have a natural homeomorphism 
between a neighborhood of $w \in X$ and a neighborhood of $w \in T$. 
The angle function is then copied from the original Hubbard Tree ${\bf H}$. 
At critical points, it is enough to extend the coordinate functions $\phi_v$ in a compatible way; 
the angle between edges can be read from this. 
We proceed as follows. 
Let $v \in T_\alpha$ be critical. 
We will define the coordinate 
$\phi_{i_\alpha(v)}$ at ${i_\alpha(v)}\in X$ as follows. 
By definition (compare $\S$II.3.6) there is a $k$ such that $\ell \in E_v$ belongs to $T_\alpha$ if and only if 
$\phi_v(\ell) \in [{k \over \delta(v)},{k+1 \over \delta(v)})$. 
Now, an edge $\ell$ incident at $i_\alpha(v)$ must belong to a unique $H_\alpha$ 
and therefore corresponds to a unique edge $\ell' \in E_{\tau(v)}$ 
in the original tree $T$. 
Define $\phi_{i_\alpha(v)}(\ell)={k + \phi_{\tau(v)}(\ell') \over \delta(v)}$. 
\smallskip
As no new periodic vertices are added the tree is still expanding. 
At periodic Julia vertices no new edges are added (compare $\S$II.3.3). 
Therefore, we have a Hubbard Tree which is homogeneous by Corollary 2.3 and satisfies the required properties. 
\medskip
To prove uniqueness, 
we note that any other local coordinate system 
$\{\phi_v\}_{v \in V}$ is also canonically present in the new tree constructed. 
It follows from Lemma 2.2 that the corresponding partition with respect to this 
coordinate is independent of the starting local coordinate system. 
\endofproof
\bigskip
{\bf 2.6 Definition.}
Let ${\bf [H]},{\bf [H']}$ be equivalence classes of abstract Hubbard Trees. 
We say that the equivalence class ${\bf [H']}$ of homogeneous abstract Hubbard Trees is the {\it inverse} of ${\bf [H]}$ 
(in symbols $inv({\bf H})={\bf H'}$), 
if ${\bf [H]}$ is the image of ${\bf [H']}$. 
\smallskip
Thus, by Propositions 2.3 and 2.5, $inv$ determines a one to one mapping from equivalence classes of abstract Hubbard Trees of degree $n>1$ to itself. 
Furthermore, in this new language Proposition 2.5 reads as follows. 
\bigskip
{\bf 2.7 Proposition.}
{\it Let ${\bf H}$ be an abstract Hubbard Tree, then $inv({\bf H})$ is a canonical extension of ${\bf H}$.} 
\endofproof
\bigskip
{\bf 2.8 Corollary.} 
{\it Let ${\bf H}$ be an abstract Hubbard Tree and $\omega$ a Fatou vertex. 
There is a canonical extension ${\bf H'}$ of ${\bf H}$ such that 
\smallskip
(a) $E_v=E'_v$ at all vertices of ${\bf H}$. 
\tskip
(b) For every $\ell \in E'_\omega$ with end points $\omega,v$, we have that $v$ is of Julia type, and 
$d_{{\bf H}'}(\tau^{\circ k}(\omega),\tau^{\circ k}(v))=1$ for all $k\ge 0$.}
\smallskip
{\bf Proof.} 
We apply first Proposition 1.2 and then take a finite number of `inverses' 
(Proposition 2.5). 
Finally we restrict to the tree generated by the original vertices. 
\endofproof
\bigskip
{\bf 2.9 Corollary.} 
{\it Let ${\bf H}$ be an abstract Hubbard Tree. 
Then ${\bf H}$ has a canonical extension in which all ends are of Julia type.}
\smallskip
{\bf Proof.} 
We apply first Proposition 1.2 and then take a finite number of `inverses' 
(Proposition 2.5). 
Finally we restrict to the tree generated by the required vertices. 
\endofproof
\insec 
\centerline{\medfont 3. Incidence.}
\insec 
{In this section we study from the dynamical point of view, 
how the number of edges incident at a Julia vertex can grow 
as we take inverses. 
This section is the abstract analogue of Section I.3.} 
\bigskip
{\bf 3.1 Definition.} 
Let ${\bf [H]}$ be an equivalence class of abstract Hubbard Trees. 
We define the incidence number $\nu_{\bf H}(v)$ at a vertex $v \in V$ 
as the number of connected components of $T-\{v\}$ in any underlying topological tree $T$. 
In the inverse trees $inv^{\circ m}({\bf [H]})$ 
we have also incidence functions $\nu_{{\bf H},-m}=\nu_{inv^{\circ m}({\bf [H]})}$ 
at the vertices of $inv^{\circ m}({\bf H})$. 
By definition $\nu_{{\bf H},0}(v) \le \nu_{{\bf H},-1}(v)$ for $v \in V$. 
Also by construction of $inv({\bf H})$, it follows that 
$\nu_{{\bf H},-1}(v)=\delta(v)\nu_{{\bf H},0}(\tau(v))$ 
for all vertices in $inv({\bf H})$. 
\bigskip
{\bf 3.2 Proposition.}
{\it Let ${\bf [H]}$ be an equivalence class of abstract Hubbard Trees.
For every periodic Julia vertex $v \in V$ and $m \ge 0$ 
we have $\nu_{H,0}(v)=\nu_{H,-m}(v)$.}
\medskip
{\bf Proof.} 
As $\delta(v')=1$, for every point $v' \in {\cal O}(v)$, 
no new edges are added around $v$ in the construction of $inv^{\circ m}({\bf [H]})$. 
(See also Lemma II.3.13.) 
\endofproof
\bigskip
{\bf 3.3 Corollary.}
{\it Let ${\bf [H]}$ be an equivalence class of abstract Hubbard Trees.
Let $v \in V$ be a Julia vertex such that $\tau^{\circ k}(v)$ is periodic. 
Then for every $m \ge k$ we have  
$\nu_{H,-k}(v)=\nu_{H,-m}(v)$.}
\endofproof
\bigskip
{\bf 3.4 Corollary.}
{\it Let ${\bf [H]}$ be an equivalence class of abstract Hubbard Trees.
There is a $k \ge 0$ such that for all 
$m \ge k$ we have $\nu_{H,-k}(v)=\nu_{H,-m}(v)$ at every Julia vertex  
$v \in V_{\bf H}$.}
\endofproof
\bigskip
We denote such numbers by $\nu_{{\bf H},-\infty}(v)$. 
\insec
\centerline{\medfont 4. Accesses and External Coordinates.}
\insec
{In this section we associate to every `access' at a Julia vertex an\break 
argument. 
This coordinate system will allow us to define extensions 
with `reasonable' properties. 
Combining these two results we prove that every extension of a Hubbard Tree is canonical.} 
\bigskip
{\bf 4.1 Definition.} (Compare Definition II.1.3.) 
Let ${\bf H}$ be an abstract Hubbard Tree. 
Given $\ell,\ell' \in E_v$ consecutive in the cyclic order, 
we say that $(v,\ell,\ell')$ is an {\it access to $v$} if $\nu_{{\bf H},0}(v)=\nu_{{\bf H},-\infty}(v)$. 
If $\nu_{{\bf H},0}(v) < \nu_{{\bf H},-\infty}(v)$ we say that $(v,\ell,\ell')$ is a {\it strict pseudoaccess to $v$} in ${\bf H}$. 
Note that at Fatou vertices there are no possible accesses. 
Clearly an access at $v$ is periodic if and only if $v$ is periodic. 
These concepts extend to equivalence classes. 
\bigskip
{\bf 4.2 Lemma.}
{\it Let ${\bf H}$ be an abstract Hubbard Tree of degree n.
Then $\tau$ induces a degree n orientation preserving covering mapping between the pseudoaccesses of the trees $inv({\bf H})$ and ${\bf H}$. 
Furthermore, accesses in $inv({\bf H})$  map to accesses in ${\bf H}$.} 
\medskip
{\bf Proof.}
If $(v,\ell,\ell')$ is a pseudoaccess in $inv({\bf H})$,  
by construction\break 
$(\tau(v),\tau_v(\ell),\tau_v(\ell'))$ is a pseudoaccess in ${\bf H}$.
Clearly this is $n$ to 1, and order preserving by construction. 
The second part is obvious. 
\endofproof
\bigskip
{\bf 4.3 Proposition.}
{\it Let ${\bf H}$ be a homogeneous abstract Hubbard Tree of degree $n>1$ with at least one Julia vertex. 
There exist an embedding $\phi_{\bf H}$ of the accesses of ${\bf H}$ into ${\bf R/Z}$ 
such that the induced map between accesses becomes multiplication by n (modulo 1). 
Furthermore $\phi_{\bf H}$ is uniquely defined up to a global addition of a  multiple of $1/(n-1)$.}
\medskip
{\bf Proof.} 
Instead of proving directly that we can assign an argument to each access of ${\bf H}$, 
we will prove this fact in a larger tree $inv^{\circ m}({\bf H})$, where $m$ is big enough. 
The result follows then by restriction (compare Lemma B.1.7 and Corollary B.2.8 in Appendix B). 
\smallskip
By Lemma 4.2 the induced map between accesses is an orientation preserving covering of degree $n$. 
In order to be able to assign an argument to each access we must prove that this map is expanding 
(compare Appendix B). 
Take two consecutive periodic accesses ${\cal A}_i=(v_i,\ell_i,\ell_i')$ in ${\bf H}$ ($i=0,1$). 
The idea is to show that for some $m$ big enough, 
these accesses are not consecutive in $inv^{\circ m}({\bf H})$. 
As no new periodic vertices are added in the construction of 
$inv^{\circ m}({\bf H})$, 
we have no new periodic accesses and the conditions of Lemma B.1.7 are trivially satisfied; 
this will establish the result. 
We distinguish between $v_0 = v_1$ and $v_0 \ne v_1$. 
\smallskip
If $v_0 = v_1$ then $\ell_0 \prec \ell_0' \preceq \ell_1 \prec \ell_1' \preceq \ell_0$ at $E_{v_0}$. 
It is enough to find an $m \ge 0$ such that $inv^{\circ m}({\bf H})$ has an access in 
the `branch' ${\cal B}_{v_0,inv^{\circ m}({\bf H})}(\ell'_0)$. 
If there is a Julia vertex in ${\cal B}_{v_0,{\bf H}}(\ell'_0)$ this is obvious by Corollary 3.4. 
If not, $\ell'_0$ has end points $v_0,\omega$ where $\omega$ is a Fatou point. 
Now the edge $\ell_0'$ corresponds to an argument 
in the coordinate $\phi_{\omega}$ at $\omega$; 
as $\omega$ eventually maps to a critical point, we can find an argument 
$\theta\ne \phi_{\omega}(\ell_0')$ which 
eventually maps to the same argument as $\phi_{\omega}(\ell)$ under successive multiplication by 
$deg_{\tau^{\circ i}(\omega)}$ modulo $1$ 
(compare diagram $(2)$ in $\S$II.3.3). 
It follows that for some $m$ big enough, 
there is an $\ell' \in E'_{\omega}$ such that 
$\phi_{\omega}(\ell')=\theta$. 
The result then follows easily from Corollary 2.8. 
(Alternatively, we can use Corollary 2.9.)
\smallskip
Now let $v_0,v_1$ be different periodic Julia points. 
By Lemma II.3.7, for some $m>0$ there is a vertex $v'$ of $inv^{\circ m}({\bf H})$ in $[v_0,v_1]_T$ for otherwise ${\bf H}$ will not be expanding.  
If $v'$ is a Julia vertex we proceed as above. 
Otherwise, we let $(v',\ell,\ell')$ be the pseudoaccess (in $inv^{\circ m}({\bf H})$ at the Fatou vertex $v'$) 
between ${\cal A}_0$, ${\cal A}_1$ in the cyclic order. 
We take an argument $\theta$ between $\phi_\omega(\ell)$ and $\phi_\omega(\ell')$ which eventually maps to the same argument as $\phi_\omega(\ell)$ and proceed as in the last paragraph. 
\endofproof
\bigskip
{\bf 4.4.} 
As every abstract Hubbard Tree ${\bf H}$ of degree $n>1$ has a canonical extension satisfying the conditions of Proposition 4.3, 
we can associate to every access a coordinate compatible with the dynamics. 
Such map $\phi_{\bf H}$ is called an {\it external coordinate}. 
In practice, this will correspond to the argument of the external ray landing throughout this access.  
\smallskip
Now let $\theta \mapsto m_n(\theta) \mapsto \dots \mapsto m_n^{\circ k}(\theta)=\theta$, 
be a periodic orbit under the standard $n$-fold multiplication in ${\bf R/Z}$. 
The question is whether there is a canonical extension of ${\bf H}$ at which accesses corresponding to the arguments  
$\{\theta,m_n(\theta),\dots,m_n^{\circ k-1}(\theta)\}$ are present. 
For this we have the following proposition. 
\bigskip
{\bf 4.5 Proposition.}
{\it Let ${\bf H}$ be a homogeneous abstract Hubbard Tree with at least one Julia vertex. 
For any election of external coordinate 
$\phi_{\bf H}$ and periodic orbit 
$\theta \mapsto m_n(\theta) \mapsto \dots \mapsto m_n^{\circ k}(\theta)=\theta$ 
under n-fold multiplication in ${\bf R/Z}$, 
there is a canonical extension of ${\bf H}$ in which accesses corresponding to 
$\{\theta,m_n(\theta),\dots,m_n^{\circ k-1}(\theta)\}$ are present.}
\medskip
{\bf Proof.}
Using Corollary 2.8 we may assume that the distance between two Fatou vertices in never equal to 1; 
and furthermore, whenever the distance between a Fatou and a Julia vertices is one, 
so is the distance between all their iterates. 
Also, because of Corollary 2.9 we may assume without loss of generality that no Fatou vertex is an end. 
We assume that there are no accesses to which 
we can associate the referred periodic 
orbit and construct a canonical extension of this tree. 
\smallskip 
{\bf Case 1.} 
The easiest way to construct extensions with periodic orbits of period $k$ 
is whenever
there is a Fatou periodic orbit of period dividing $k$. 
Suppose the total degree of such critical cycle is $d$.   
In this case, for all arguments of period $k$ under $m_d$ we can include 
an edge which correspond in local coordinates to this argument and a periodic vertex 
(if they are not already present). 
When this is done simultaneously at all Fatou vertices of the cycle we clearly get a 
new expanding Hubbard Tree. 
Clearly this construction is canonical. 
If the required accesses are present in this canonical extension, we are done; 
otherwise we have to work harder. 
\smallskip
To continue the general case, 
first note that Corollary 3.4 guarantees that for $m$ big enough 
$\nu_{-m,{\bf H}}(v)=\nu_{-\infty,{\bf H}}(v)$ 
at every original vertex $v \in V_{\bf H}$. 
We will only keep track of the following information: 
the original tree ${\bf H}$ and all these accesses of $inv^{\circ m}({\bf H})$ 
at vertices $v \in V_{\bf H}$ in the original tree
(we have `pruned' the tree $inv^{\circ m}({\bf H})$). 
In this case if $\ell \in E^m_v$ but $\ell \not\in E_v$ 
(i.e, if the germ $\ell$ at $v$ in the tree $inv^{\circ m}({\bf H})$ 
is not present in ${\bf H}$) 
we say that the tree $inv^{\circ m}({\bf H})$ was {\it pruned} at $\ell$. 
\smallskip
\smallskip
Let $\{\gamma_1,\dots,\gamma_\alpha\}$ be the arguments of all such accesses ordered counterclockwise. 
Working if necessary in a canonical extension, 
we may further suppose that the Lebesgue measure of $(\gamma_i,\gamma_{i+1})$ 
is at most $1/n^{2k+2}$. 
(In fact, we may work in an inverse $inv^{\circ l}({\bf H})$ with $l$ big enough 
thanks to the expansiveness of $m_n$ in ${\bf R/Z}$.) 
It follows that $(\gamma_i,\gamma_{i+1})$ 
contains at most one periodic orbit of period diving $2k$ in its closure. 
In particular, each $m_n^{\circ i}(\theta)$ belongs to an interval 
$(\theta^+_i,\theta^-_i)$ with $\theta^\pm_i$ strictly preperiodic. 
It follows that the vertices $v^+_{\theta_i},v^-_{\theta_i}$ 
at the respective accesses are not periodic. 
\smallskip 
Suppose first that $v^+_{\theta_0} \ne v^-_{\theta_0}$. 
By further subdividing the tree 
(for example by taking an extra $k$ inverses and 
restricting to the vertices in the original underlying topological tree), 
we may suppose that for any edge $\ell$, 
the iterated maps 
$\tau^{\circ i}|_{\ell}$ are one to one ($i=1, \dots k$). 
\smallskip
{\bf Case 2.} 
Suppose $[v^+_{\theta_0},v^-_{\theta_0}] \subset \tau^{\circ k}([v^+_{\theta_0},v^-_{\theta_0}])$. 
It follows from standard techniques for subshifts of finite type 
that we can canonically extend the vertices of the tree so that it includes an orbit of period $k$ or $k/2$ with 
$v_{m_n^{\circ i}(\theta)} \in [v^+_{\theta_i},v^-_{\theta_i}]$. 
Because $v^\pm_{\theta_i}$ are strictly preperiodic, 
the expansive condition for the new set of vertices is trivially satisfied . 
Therefore any access at $v_{m_n^{\circ i}(\theta)}$ belonging to the set 
$(\theta^+_i,\theta^-_i)$ 
should have an associated argument of period dividing $2k$. 
By construction this argument can only be $m_n^{\circ i}(\theta)$. 
\smallskip
{\bf Case 3.} 
Suppose $[v^+_{\theta_0},v^-_{\theta_0}] \cap \tau^{\circ k}([v^+_{\theta_0},v^-_{\theta_0}])=[v_1,v_2]$. 
Then the vertices $v_1,v_2$ belong to the interval $[v^+_{\theta_0},v^-_{\theta_0}]$. 
Now, by hypothesis this last interval contains no vertex of Julia type 
(for otherwise after completing the accesses at such vertex, we will have that 
$\theta^+_0$ and $\theta^-_0$ are not consecutive in the cyclic order) 
and at most one vertex $w$ of Fatou type. 
It follows that $v_1=w$ and that 
$v_2$ equals either $v^+_{\theta_0}$ or $v^-_{\theta_0}$. 
In either case we get 
$d(\tau^{\circ k}(v^+_{\theta_0}),\tau^{\circ k}(v^-_{\theta_0})) \ge 3$. 
However, by assumption this is impossible since 
$d(\tau^{\circ k}(v^+_{\theta_0}),\tau^{\circ k}(w))=d(\tau^{\circ k}(v^-_{\theta_0}),\tau^{\circ k}(w))=1$ implies 
$d(\tau^{\circ k}(v^+_{\theta_0}),\tau^{\circ k}(v^-_{\theta_0})) \le 2$. 
\smallskip
{\bf Case 4.} 
Suppose $[v^+_{\theta_0},v^-_{\theta_0}]$ intersects 
$\tau^{\circ k}[v^+_{\theta_0},v^-_{\theta_0}]=
[\tau^{\circ k}(v^+_{\theta_0}),\tau^{\circ k}(v^-_{\theta_0})]$ 
at an interior vertex $w \in [v^+_{\theta_0},v^-_{\theta_0}]$. 
It follows from the preliminary discussion in case 3 that $w$ is a Fatou vertex. 
This Fatou vertex should be periodic of period dividing $k$ 
because otherwise 
$\tau^{\circ k}(w) \ne w$ belongs to 
$[\tau^{\circ k}(v^+_{\theta_0}),\tau^{\circ k}(v^-_{\theta_0})]$ 
and therefore 
$d(\tau^{\circ k}(v^+_{\theta_0}),\tau^{\circ k}(v^-_{\theta_0})) \ge 3$,  
which can be shown to be impossible as in case 3. 
\smallskip
Denote by $\ell^\pm_k$ the edges $[w,\tau^{\circ k}(v^\pm_{\theta_0})]$ with local coordinates $\alpha_k^\pm$ at $w$, and  
by $\ell^\pm_0$ the edges $[w,v^\pm_{\theta_0}]$ with local coordinates $\alpha_0^\pm$. 
Clearly $(\alpha^+_0,\alpha^-_0) \subset (\alpha^+_k,\alpha^-_k)$ because this is the only ordering compatible with the order of the accesses. 
Denote by $d$ the local degree of $\tau^k$ at $w$. 
\smallskip
{\bf Claim.} 
{\it $m_d$ maps $(\alpha^+_0,\alpha^-_0)$ homeomorphicaly onto 
$(\alpha^+_k,\alpha^-_k)$.} 
\smallskip
In fact, if this is not the case, 
in some inverse tree there is an edge $\ell'=[w,v']$ 
with corresponding  argument $\phi_w(\ell') \in  (\alpha^+_0,\alpha^-_0)$ and with 
$\tau^{\circ k}(\ell')=\tau^{\circ k}(\ell^+_0)$. 
It follows that after completing the access at the vertex $v'$ there is an access 
with corresponding argument $\beta \in (\theta^+_0,\theta^-_0)$ such that 
$m_n(\beta)=m_n(\theta_{\theta^+_0})$. 
But this implies that the interval $(\theta^+_0,\theta^-_0)$ has Lebesgue measure at least $1/n^k$, 
which is a contradiction. 
\smallskip
To finish the proof of case 4, we notice that the claim implies that 
$m_d$ has a fixed point inside the interval $(\alpha^+_0,\alpha^-_0)$.  
Therefore we are in case 1. 
\medskip
{\bf Case 5.} 
Suppose the intervals $[v^+_{\theta_0},v^-_{\theta_0}]$ and 
$[\tau^{\circ k}(v^+_{\theta_0}),\tau^{\circ k}(v^-_{\theta_0})]$ 
have disjoint interiors. 
In this case we consider the subtree generated by the vertices 
$v^\pm_0$ and $\tau^{\circ k}(v^\pm_0)$ to notice that there is vertex $v$ strictly contained in the interior of 
$[\tau^{\circ k}(v^+_{\theta_0}),\tau^{\circ k}(v^-_{\theta_0})]$. 
Also there is an edge $\ell$ at this vertex 
such that $v^\pm_{\theta_0} \in {\cal B}(\ell)$ the branch at $\ell$. 
In fact, this follows from the ordering of accesses. 
This implies in particular that for some inverse of the tree 
there is a vertex $v' \in [v^+_{\theta_0},v^-_{\theta_0}]$ 
with $\tau^{\circ k}(v')=v$. 
Also, we can find an edge $\ell'$ at $v'$ which 
maps locally to $\ell$ under $\tau^k$. 
If $v'$ is of Julia type, 
there are consecutive accesses (after completing the accesses) 
at $v'$ with associated arguments $\theta_{\cal A}$ and $\theta_{\cal B}$ 
such that $\theta \in (\theta_{\cal A},\theta_{\cal B}) \subset (\theta^+_0,\theta^-_0)$. 
If $v'$ is of Fatou type, 
there is a Julia vertex $v_1$ in the branch ${\cal B}(\ell')$ such that 
(after restriction to the a tree which only includes this vertex in such branch) 
there are two consecutive accesses with that property described above. 
In fact, 
these two properties follow immediately from the fact that accesses at $v'$ (respectively at $v_1$) map to accesses at $\tau^{\circ k}(v')$ (respectively at $\tau^{\circ k}(v_1)$), 
and that $(\theta^+_0,\theta^-_0)$ has Lebesgue measure at most $1/n^{2k+2}$. 
\smallskip
In either case we have reduced the problem to case 6. 
\medskip
{\bf Case 6.} 
Suppose now that $v^+_{\theta_0} = v^-_{\theta_0}$. 
After taking inverses and restricting if necessary we may suppose that 
$\tau^{\circ i}(v^\pm_{\theta_0})=v^+_{\theta_i}$ for $i=0, \dots k-1$. 
Thus, the accesses ${\cal A}^+_i$ and ${\cal A}^-_i$ 
with external arguments $\theta^+_i,\theta^-_i$ 
share an edge $\ell_i$. 
As there is no further access with argument in $(\theta^+_i,\theta^-_i)$ 
it follows that some tree $inv^{\circ m}({\bf H})$ was ``pruned" at $\ell_i$. 
In this way, the required extension is achieved by adding the vertices 
$v_{m_n^{\circ i}}(\theta)$ 
at the other end of $\ell_i$. 
Note that the extension is canonical because 
for any extension including the vertex $v^+_{\theta_i}$, the vertex 
$v_{m_n^{\circ i}}(\theta)$ should belong to the branch $\ell_i$, 
and thus, according to Lemma II.3.13  
these periodic vertices should be ends.  
\endofproof
\bigskip
{\bf 4.6 Corollary.} 
{\it Every extension of an abstract Hubbard Tree is\break 
canonical.}
\medskip
{\bf Proof.}
Given any extension we assign to every periodic access its canonical argument 
(compare Proposition 4.3). 
Then starting with the minimal tree
we add all these periodic orbits according to Proposition 4.5. 
Finally, we take a finite number of inverses and restrict if necessary. 
\endofproof
\smallskip
\insec 
\centerline{\medfont 5. From Hubbard Trees to Formal Critical Portraits.}
\insec
{Using canonical extensions we will mimic the constructions of critical portraits from the first part of this work. 
For the main defenitions and results see Appendix A} 
\bigskip
{\bf 5.1 Extending the tree.} 
Let ${\bf H}$ be an abstract Hubbard Tree of degree $n>1$. 
We start with a canonical extension ${\bf H'}$ of ${\bf H}$ as in Corollary 2.8; 
i.e, we require from this extension that 
if $\omega$ is a Fatou point and $\ell \in E_{\omega}$, 
then for the endpoints $\omega,v$ of $\ell$ we must have that 
$v$ is a Julia vertex, and 
$d_{\bf H'}(\tau^{\circ k}(v),\tau^{\circ k}(\omega))=1$ for all $k \ge 0$. 
\smallskip
We fix local coordinates $\{\phi_v\}_{v \in V}$. 
For any critical cycle we extend the tree by adding an edge and a vertex at every 0 argument (if they are not present).  
Next, for any Fatou vertex $\omega$ we proceed as follows. 
Inductively suppose that the 0 edge is present in the local coordinate of $\tau(\omega)$. 
We insert a new vertex and edge (if they are not present) at every argument of 
$\phi^{-1}_{\tau(\omega)}(0)$. 
Then we use Corollary 3.4 to guarantee that pseudoaccesses defined at such points are indeed accesses. 
We call any extension satisfying the above conditions {\it supporting} 
(compare $\S$I.2).
\medskip
Let $\omega$ be a Fatou vertex, 
an access $(v,\ell',\ell)$  is said to support $\omega$ if 
$\ell$ has endpoints $\omega,v$ and 
$d_{\bf H}(\tau^{\circ k}(\omega),\tau^{\circ k}(v))=1$ for all $k \ge 0$. 
Clearly $\tau(v,\ell',\ell)=(\tau(v),\tau_v(\ell'),\tau_v(\ell))$ supports $\tau(\omega)$. 
An access $(v,\ell',\ell)$ which supports the Fatou critical point $\omega$ 
will be denoted by ${\cal D}(\omega,\ell)$
\bigskip
{\bf 5.2 Constructing marked accesses.} 
Let ${\bf H}$ be a supporting abstract Hubbard Tree. 
Using Corollary 3.4 we pick an inverse $inv^{\circ m}({\bf H})$ such that 
at every $v \in V$ we have $\nu_{{\bf H},-m}(v)=\nu_{{\bf H},-\infty}(v)$. 
From this it is easy to chose hierarchic accesses as in $\S$I.2: 
\smallskip
For each critical vertex $\omega \in \Omega({\bf H})$ set 
$$
\Lambda_{\omega}=\{\ell \in E_\omega: \delta(\omega)\phi_v(\ell)=0\}
$$
(in this case the hierarchic selection is reflected in the choice of a 0 argument in the local coordinate). 
Let $\Omega({\cal F})=\{\omega^{\cal F}_{1},\dots,\omega^{\cal F}_{l}\}$ be the set of Fatou critical vertices, and  
$\Omega({\cal J})=\{\omega^{\cal J}_{1},\dots,\omega^{\cal J}_{k}\}$ the set of Julia critical vertices. 
For each $\omega \in \Omega({\cal F})$ we construct $\delta(\omega)$ {\it marked} supporting accesses to $\omega$ in the following way. 
Take $\ell \in \Lambda_\omega$ with end points $v_\ell,\omega$; 
then there is a supporting access to $\omega$ at $v_\ell$ of the form ${\cal D}(\omega,\ell)=(v_\ell,\ell',\ell)$. 
The set of such accesses for all possible $\ell \in \Lambda_\omega$ is by definition 
${\cal F}_\omega$. 
\smallskip
For each $\omega \in \Omega({\cal J})$ we construct $\delta(\omega)$ {\it marked} accesses in the following way. 
Take $\ell \in \Lambda_\omega$, 
then there is an accesses at $\omega$ of the form ${\cal E}(\omega,\ell)=(\omega,\ell,\ell')$. 
The set of such accesses for all possible $\ell \in \Lambda_\omega$ is by definition 
${\cal J}_\omega$. 
\smallskip
Note the slight difference in the construction, 
at a Julia critical vertex $v$, 
the marked accesses are at $v$. 
While for Fatou critical vertices the accesses are taken at the other end of each edge. 
\medskip
In this way we have constructed two families 
$$
\matrix{ 
{\cal F}&=\{{\cal F}_{\omega_1},\dots,{\cal F}_{\omega_l}\} \cr
{\cal J}&=\{{\cal J}_{\omega_1},\dots,{\cal J}_{\omega_k}\} \cr}
$$
of accesses. 
As these accesses correspond in the external coordinate $\phi_{\bf H}$ to arguments, 
we will not distinguish between the accesses and their corresponding argument. 
In this way we have the following (see $\S$I.3). 
\bigskip
{\bf 5.3 Proposition.} 
{\it The marking $({\cal F},{\cal J})$ is a formal critical portrait.} 
\bigskip
{\bf Proof.} 
This follows directly from the construction.
\endofproof
\bigskip
There are several trivial consequences of this construction that we want to point out. 
To simplify notation, 
the vertex at which an access ${\cal C}$ is defined 
will be denoted by $v_{\cal C}$. 
The proof in all cases is the same: 
by removing the edge $\ell$ we are left with two connected pieces. 
\bigskip
{\bf 5.4 Lemma.}
{\it Let $\omega$ be a Fatou critical vertex. 
If $v_{\cal C} \in {\cal B}_{{\bf H},\omega}(\ell)$, then 
for all  
$\ell' \in \Lambda_\omega-\{\ell\}$ we have
${\cal D}(\omega,\ell') \prec {\cal C} \preceq {\cal D}(\omega,\ell)$.}
\endofproof
\bigskip
{\bf 5.5 Lemma.}
{\it Let $\omega$ be a Julia critical vertex, 
and ${\cal C}$ an access at\break 
$v_{\cal C} \in {\cal B}_{{\bf H},\omega}(\ell)-\{\omega\}$. 
Then for any accesses ${\cal A},{\cal A}'$ at $\omega$ we have either\break 
${\cal A} \prec {\cal C} \prec {\cal A}'$ or 
${\cal A}' \prec {\cal C} \prec {\cal A}$.} 
\endofproof
\bigskip
{\bf 5.6 Lemma.}
{\it Suppose $\omega$ is a Fatou critical vertex and  
let $\ell \not\in \Lambda_\omega$. 
If ${\cal C}$ an access at 
$v_{\cal C} \in {\cal B}_{{\bf H},\omega}(\ell)$, 
then for any $\ell',\ell'' \in \Lambda_\omega$ we have either 
${\cal D}(\omega,\ell') \prec {\cal C} \prec {\cal D}(\omega,\ell'')$ or 
${\cal D}(\omega,\ell'') \prec {\cal C} \prec {\cal D}(\omega,\ell')$.} 
\endofproof
\insec
\centerline{\medfont 6. From Hubbard Trees to Admissible Critical Portraits.}
\insec 
{In this section we prove that the formal critical portrait constructed above 
is also admissible.  
For this we must verify conditions $(c.6),(c.7)$ in $\S$A.2.7. 
We first verify condition $(c.6)$. 
The verification of condition $(c.7)$, 
will also show that any polynomial with critical marking $({\cal F},{\cal J})$ has Hubbard Tree equivalent to this starting one. 
In this way the main Theorem {\bf A} will follow.} 
\bigskip
{\bf 6.1 Proposition.}
{\it The formal critical portrait $({\cal F},{\cal J})$ is an admissible critical portrait.}
\medskip
{\bf Proof.} 
This follows from Corollaries 6.4 and 6.9 below.
\endofproof
\bigskip
{\bf 6.2 Lemma.} 
{\it Let ${\cal A}_i,{\cal B}_i$ be accesses at $v_i$ for $i=1,2$ with $v_1 \ne v_2$. 
Then $\{{\cal A}_1,{\cal B}_1\}$, and $\{{\cal A}_2,{\cal B}_2\}$ are unlinked.}
\medskip
{\bf Proof.} 
This follows from the fact that $\{{\cal A}_2,{\cal B}_2\}$ are defined in the same connected component of $T-\{v_1\}$. 
\endofproof
\bigskip
{\bf 6.3 Lemma.} 
{\it Let ${\cal A},{\cal A}'$ be periodic accesses. 
If either $S^+({\cal A})=S^+({\cal A}')$ or  
$S^-({\cal A})=S^-({\cal A}')$, 
then $v_{\cal A}=v_{{\cal A}'}$}.
\medskip
{\bf Proof.} 
By contradiction suppose $v_{\cal A} \ne v_{{\cal A}'}$. 
We distinguish two cases. 
\smallskip
Suppose $\tau^{\circ k}|_{[v_{\cal A},v_{{\cal A}'}]_T}$ is injective for all $k \ge 1$. 
In this case there is a periodic Fatou vertex $v \in [v_{\cal A},v_{{\cal A}'}]_T$, 
because otherwise the tree will not be expanding. 
Let $d>1$ be the degree of the critical cycle $v_0=v \to v_1 \dots \to v_n=v_0$. 
There are exactly two different edges $\ell,\ell' \in E_v$ contained in $[v_{\cal A},v_{{\cal A}'}]_T$. 
The dynamics of these edges must be periodic by Lemma II.3.7. 
We write $\phi_v(\ell),\phi_v(\ell')$ in base $d$ expansion. 
As they are not equal by hypothesis, we may suppose that the first coefficient in the expansions are different. 
As $d$ is the product of the degrees of the vertices in the cycle, 
we may suppose then that when multiplying by $\delta(v_0)$ they have different integer part. 
But in this way by Lemma 5.6 we will have 
$\pi_0(S^+({\cal A})) \ne  \pi_0(S^+({\cal A}'))$. 
(In fact, for $\epsilon>0$ small enough, 
the arguments $\phi_{\bf H}({\cal A})$ and $\phi_{\bf H}({\cal A}')$ 
belong to different connected components of 
${\bf R/Z}-\{\phi_{\bf H}({\cal D}(v,\ell)): \ell \in \Lambda_v\} =\break 
{\bf R/Z}-{\cal F}_v$.) 
But implies that  
$S^+({\cal A}) \ne  S^+({\cal A}')$. 
If we consider instead of $\phi_v$ the `coordinate' $1-\phi_v$ 
the same reasoning give us  
$S^-({\cal A}) \ne  S^-({\cal A}')$.
\smallskip
Suppose now that $\tau|_{[v_{\cal A},v_{{\cal A}'}]_T}$ is not locally one to one near $\omega$. 
If $\omega$ is a Julia critical vertex the result follows from Lemma 5.5. 
If $\omega$ is a Fatou critical vertex, by Lemma 5.6 we always have 
$\pi_0(S^-({\cal A})) \ne  \pi_0(S^-({\cal A}'))$ and thus 
$S^-({\cal A}) \ne  S^-({\cal A}')$.
\smallskip
If neither ${\cal A}$ nor ${\cal A}'$ support $\omega$, 
again by Lemma 5.6 
$\pi_0(S^+({\cal A})) \ne  \pi_0(S^+({\cal A}'))$. 
We start though by assuming that ${\cal A}$ is a marked access associated with $\omega$. 
By Hypothesis there is a preperiodic marked access ${\cal C} \in {\cal F}_\omega$ 
(and therefore such that $\tau({\cal C})=\tau({\cal A})$) 
with   
$v_{\cal C} \in [v_{\cal A},v_{{\cal A}'}]_T$. 
Thus $\tau^{\circ k}|_{[v_{\cal C},v_{{\cal A}'}]_T}$ eventually maps into $[v_{\cal A},v_{{\cal A}'}]_T$. 
It follows there is a point $\omega' \in [v_{\cal C},v_{{\cal A}'}]_T$ that eventually maps to $\omega$. 
Working if necessary in a canonical extension $inv^{\circ k}({\bf H})$ 
we may assume without loss of generality 
that $\omega' \in V$. 
But then by Lemma II.3.7 for some $i \ge k$, 
$\tau^{\circ i}|_{[\omega,\omega']_T}$ is not locally one to one near some point $\omega''$. 
If $i$ is minimal, 
neither of the periodic accesses $\tau^{\circ i}({\cal A})=\tau^{\circ i}({\cal C})$ nor $\tau^{\circ i}({\cal A}')$ 
can support the critical point $\tau^{\circ i-1}(\omega'')$ if it is of Fatou type. 
It follows from the previous reasoning that 
$S^+(\tau^{\circ i-1}({\cal A})) \ne S^+(\tau^{\circ i-1}({\cal A}'))$, 
and therefore 
$S^+({\cal A}) \ne S^+({\cal A}')$. 
\endofproof
\bigskip
{\bf 6.4 Corollary.}
{\it The formal critical portrait $({\cal F},{\cal J})$ satisfies condition (c.6).}
\medskip
{\bf Proof.} Let ${\cal A}$ be a periodic marked access. 
Suppose there is a periodic argument 
$\lambda$ such that $S^+(\lambda)=S^+({\cal A})$. 
By Proposition 4.5 
we can assume that there is an accesses corresponding to $\lambda$. 
By Lemma 6.3 this access is supported at $v_{\cal A}$. 
By Lemma 5.4 this access can only be ${\cal A}$.
\endofproof
\medskip
{\bf 6.5 Lemma.} 
{\it Let $v_{\cal A}=v_{{\cal A}'}$ be a non critical Julia vertex. 
Then ${\cal A}$ and ${\cal A}'$ have the same left address, i.e, 
$\pi_0(S^-({\cal A}))=\pi_0(S^-({\cal A}'))$.}
\smallskip
{\bf Proof.}
If ${\cal E,E}'$ are marked accesses associated to the same Julia critical vertex, 
Lemma 6.2 implies that $\{{\cal A,A}'\}$,$\{{\cal E,E}'\}$ are unlinked. 
\smallskip
If ${\cal D,D}'$ are marked accesses associated to the same Fatou critical vertex, 
we distinguish if $v_{\cal A}$ equals $v_{\cal D}$ or not. 
If $v_{\cal A} \ne v_{\cal D},v_{{\cal D}'}$ then 
clearly $\{{\cal A,A}'\}$,$\{{\cal D,D}'\}$ are unlinked because 
the regulated path $[v_{\cal D},v_{{\cal D}'}]_T$ does not contain $v_{\cal A}$. 
If $v_{\cal A} = v_{\cal D}$ then by Lemma 5.5 
${\cal D}' \prec {\cal A} \prec {\cal A}' \preceq {\cal D}$. 
\smallskip
All these facts together mean by definition that the accesses 
${\cal A}$ and ${\cal A}'$ have the same left address, i.e, 
$\pi_0(S^-({\cal A}))=\pi_0(S^-({\cal A}'))$. 
\endofproof
\medskip
{\bf 6.6 Lemma.} 
{\it Let ${\cal B}$ be an access at a Julia critical vertex v. 
Then there is a marked access ${\cal E}$ at v, 
such that $\pi_0(S^-({\cal E}))=\pi_0(S^-({\cal B}))$.}
\smallskip
{\bf Proof.} 
Take consecutive ${\cal E}$, ${\cal E}'$ marked accesses at $v$, such that\break 
${\cal A}' \prec {\cal E} \preceq {\cal A}$. 
Using Lemma 6.2 and the same reasoning as in Lemma 6.5 
we get $\pi_0(S^-({\cal A}))=\pi_0(S^-({\cal E}))$.
\endofproof
\medskip
{\bf 6.7 Corollary.}
{\it Suppose $\pi_0(S^-({\cal A}))=\pi_0(S^-({\cal A}'))$. 
Then $v_{\cal A}=v_{{\cal A}'}$ 
if and only if 
$v_{\tau({\cal A})}=v_{\tau({\cal A}')}$.}
\smallskip
{\bf Proof.} 
One direction is obvious. 
On the other hand, we may assume that $v_{\tau({\cal A})}$ has $n$ inverses in the tree counting multiplicity. 
As there are only $n$ possible choices of addresses, 
the result follows combining Lemmas 6.3, 6.5, 6.6. 
\endofproof
\smallskip
{\bf 6.8 Proposition.}
{\it $v_{\cal A}=v_{{\cal A}'}$ if and only if 
$S^-({\cal A}) \sim_l S^-({\cal A}')$.}
\smallskip
{\bf Proof.} 
First suppose $S^-({\cal A})\sim_l S^-({\cal A}')$. 
It is enough to prove that if $S^-({\cal A})\approx S^-({\cal A}')$ 
then $v_{\cal A}=v_{{\cal A}'}$. 
If $S^-({\cal A}) = S^-({\cal A}')$ this follows from Lemma 6.3 and Corollary 6.7. 
In the other case the result follows from this fact, 
Lemma 6.6 and again Corollary 6.7.
\smallskip
Suppose now $v_{\cal A}=v_{{\cal A}'}$. 
Let $m \ge 0$ be the smallest integer such that 
$\tau^{\circ m}(v_{\cal A})$ does not contain in its forward orbit a critical vertex. 
The proof will be in induction in $m$. 
For $m=0$ this is Lemma 6.5. 
Suppose now that the result holds for $m-1$. 
This implies that all accesses at
 $\tau(v_{\cal A})$ have equivalent symbol sequences. 
If $v$ is not critical we use again Lemma 6.5. 
If $v$ is critical we use Lemma 6.6. 
\endofproof
\smallskip
{\bf 6.9 Corollary.}
{\it The formal critical portrait $({\cal F},{\cal J})$ satisfies condition (c.7).}
\endofproof
\eject
\insec 
\centerline{\medfont 7. Proof of the Theorem A.}
\insec 
The admissible critical portrait $({\cal F},{\cal J})$ determines a unique (up to affine conjugation) polynomial $P$ 
with marking $(P,{\cal F},{\cal J})$ by Theorem A.2.9. 
By Propositions 6.8 and A.2.12 its Hubbard Tree is the starting one. 
The angle function at Fatou vertices are the starting ones because of Proposition 2.7, and Corollaries 2.8 and B.2.5
\endofproof
\bigskip

%% file: apen.ht.tex
\inchap 
\centerline{\bigfont Appendix A} 
\centerline{\bigfont Critical Portraits.}
\inchap
{\sl We follow closely the exposition in [P1] about critical portraits for postcritically finite polynomials. 
Proofs of all statements and further details can be found in said work.} 
\insec 
\centerline{\medfont $\S$A.1 Construction of Critically Marked Polynomials.} 
\insec
{\bf A.1.1 Supporting arguments.} 
Let $P$ be a $PCF$ Polynomial. 
Given a Fatou component $U$ and a point $p \in \partial U$, 
there are only a finite number of external rays $R_{\theta_1},\dots,R_{\theta_k}$ landing at $p$. 
These rays divide the plane in $k$ regions.
We order the arguments of these rays in counterclockwise cyclic order 
$\{\theta_1,\dots,\theta_k\}$, 
so that $U$ belongs to the region determined by 
$R_{\theta_1}$ and $R_{\theta_2}$ ($\theta_1=\theta_2$ if a single ray lands at $p$). 
The argument $\theta_1$ (respectively the ray $R_{\theta_1}$) is by definition the 
{\it (left) supporting argument (respectively the (left) supporting ray) of the Fatou component U}. 
In a completely analogous way we can define right supporting rays. 
Note that an argument supports at most one Fatou component. 
Furthermore, by definition, 
given a Fatou component $U$, 
at every boundary point $p$ lands a supporting ray for $U$. 
\smallskip
{\bf Definition.} 
Given an external ray $R_\theta$ supporting the Fatou component $U(z)$ with center $z$, 
we extend $R_\theta$ by joining its landing point with $z$ by an internal ray, 
and call this set an {\it extended ray $\hat R_\theta$ with argument $\theta$}. 
\bigskip
Given a postcritically finite polynomial $P$, which we assume to be monic and centered, 
we associate to every critical point a finite subset of ${\bf Q/Z}$ and construct a 
{\it critically marked polynomial} 
$(P,{\bf \cal F}=\{{\cal F}_1,\dots,{\cal F}_{n_F}\},{\bf \cal J}=\{{\cal J}_1,\dots,{\cal J}_{n_J}\})$. 
Here ${\cal F}_k$ would be the set of arguments associated with the critical point $z^F_k$ in the Fatou set, 
and ${\cal J}_k$ would be the set associated with the critical point $z^J_k$ in the Julia set. 
We remark that given a polynomial its critical marking is not necessarily unique. 
Also note that one of these two families may be empty if there are no critical points in the Fatou or Julia sets.
In the following definition we will always work with left supporting rays. 
We remark that we could equally well work with the right analogue, 
but there must be the same choice throughout. 
Also, multiplication by $d$ modulo $1$ in {\bf R/Z} will be denoted by $m_d$. 
\bigskip
\bigskip
{\bf A.1.2 Construction of ${\cal F}_k$.} 
First we consider the case in which a given Fatou critical point $z=z^F_k$ is periodic. 
Let $z=z^F_k \mapsto P(z) \mapsto ... \mapsto P^{\circ n}(z)=z$ be a critical cycle of period $n$ and degree ${\cal D}_z>1$ 
(by definition the degree of a cycle is the product of the local degree of all elements in said cycle). 
We construct the associated set ${\cal F}_k$ for every critical point in the cycle simultaneously. 
Denote by $d_z$ be the local degree of $P$ at $z$. 
We pick any periodic point $p_z \in \partial U(z)$ of period dividing $n$ 
(which is not critical, because is periodic and belongs to the Julia set $J(P)$) 
and consider the supporting ray $R_\theta$ for this component $U(z)$ at $p_z$. 
Note that this choice naturally determines a periodic supporting ray for every Fatou component in the cycle. 
The period of this ray is exactly $n$. 
Given this periodic supporting ray $R_\theta$, 
we consider the $d_z$ supporting rays for this same component $U(z)$ that are inverse images of $P(R_\theta)=R_{m_d(\theta)}$. 
The set of arguments of these rays is defined to be ${\cal F}_k$. 
Keeping in mind that a preferred periodic supporting ray has been already chosen, 
we repeat the same construction for all critical points in this cycle. 
Note that as the cycle has critical degree ${\cal D}_z$, 
we can produce ${\cal D}_z-1$ different possible choices for ${\cal F}_k$. 
If ${\cal F}_k$ is the set associated with the periodic critical point $z_k$, 
there is only one periodic argument in ${\cal F}_k$ 
(namely $\theta$ as above), 
we call this angle {\it the preferred supporting argument associated with $z^F_k$.} 
Note that by definition, the period of $z_k^F$ 
equals the period of the associated preferred periodic argument. 
\smallskip
Otherwise, if $z=z^F_k$ of degree $d_z>1$, 
is a non periodic critical point in the Fatou set $F(P)$, 
there exists a minimal $n>0$ for which $w=P^{\circ n}(z)$ is critical. 
If $w$ has associated a preferred supporting ray $R_\theta$ 
(at the beginning only periodic critical points do), 
then in $P^{-n}(R_\theta)$ there are exactly $d_z$ rays that support this Fatou component $U(z)$. 
The set of arguments of these rays is defined to be ${\cal F}_k$. 
We pick any of those and call it the {\it preferred supporting argument associated with $z$}. 
We continue this process for all Fatou critical points.
\bigskip
{\bf A.1.3 Construction of ${\cal J}_k$.} 
Given $z=z_k^J$ (a critical point in $J(P)$) of degree $d_k>1$, 
we distinguish two cases. 
If the forward orbit of $z$ contains no other critical point, 
we have that for some $\theta$ (usually non unique) $R_\theta$ lands at $P(z)$. 
Now $P^{-1}(R_\theta)$ consists of $d$ different rays, 
among them exactly $d_k$ land at $z$. 
Define ${\cal J}_k$ as the set of arguments of these rays, 
and choose a {\it preferred ray}. 
Otherwise, $z$ will map in $n \ge 1$ iterations to a critical point, 
which we assume to have associated a preferred ray $R_{\theta}$. 
In the $n^{th}$ inverse image $P^{-n}(R_\theta)$ of this preferred ray, 
there are $d_k$ rays which land at $z$. 
The set of arguments of these rays is defined to be ${\cal J}_k$. 
Again we pick one of those to be preferred, 
and continue until every critical point has an associated set. 
\bigskip 
The critical marking itself gives information about 
how many iterates are needed for a given critical point to become periodic. 
For example, by construction we have the following lemma. 
\smallskip
{\bf A.1.4 Lemma.} 
{\it Let $\gamma$ be a preferred supporting argument in ${\cal F}_k$ (respectively in ${\cal J}_k$). 
Then the multiple 
$m^{\circ n}_d(\gamma)$ (with $n \ge 1)$ is periodic but $m_d^{\circ n-1}(\gamma)$ is not 
if and only if 
$z^F_k$ (respectively $z^J_k$) 
falls in exactly $n$ iterations into a periodic orbit.} 
\bigskip
{\bf Remark.} 
Note that the construction of associated sets was done in several steps. 
We first complete the choice for all critical cycles, 
and then proceed backwards. 
In both the Fatou and Julia set cases we will have to make decisions at several stages of the construction. 
Such decisions will affect the choice of the marking for all critical points found in the backward orbit of these starting ones. 
Each time that this kind of construction is made, 
we will informally say that it is a {\it hierarchic selection}. 
\insec
\centerline{\medfont $\S$A.2 The Combinatorics of Critically Marked Polynomials.} 
\insec
In order to analyze which conditions the families $({\bf \cal F},{\bf \cal J})$ satisfy, 
it is convenient to introduce some combinatorial notation. 
\bigskip
{\bf A.2.1 Definitions.} 
We say that a subset ${\Lambda} \subset {\bf R/Z}$ is a {\it (d-)preargument set} 
if $m_d({\Lambda})$ is a singleton. 
For technical reasons we assume always that $\Lambda$ contains at least two elements. 
If all elements of $\Lambda$ are rational, we say that $\Lambda$ is a {\it rational preargument set}. 
It follows by construction that whenever $(P,{\bf {\cal F},{\cal J}})$ is a marked polynomial, 
all the sets ${\cal J}_k$, and ${\cal F}_l$ are rational $d$-preargument sets.
\medskip
Consider now a family ${\bf \Lambda}=\{\Lambda_1,\dots,\Lambda_n\}$ of finite subsets of the unit circle ${\bf R/Z}$. 
The family ${\bf \Lambda}$ determines the {\it family union} set ${\bf \Lambda}^\cup =\bigcup \Lambda_i$. 
We say that any $\lambda \in {\bf \Lambda}^\cup$ is an {\it element of the family ${\bf \Lambda}$}. 
Furthermore, we can say that it is a periodic or preperiodic element of the family if it is so with respect to $m_d$. 
The set of all periodic elements in the family union will be denoted by 
${\bf \Lambda^\cup_{per}}$.
\bigskip
{\bf A.2.2 Hierarchic Families.} 
We say that a family ${\bf \Lambda}$ is {\it hierarchic} if 
for any elements in the family $\lambda,\lambda' \in {\bf \Lambda}^\cup$, 
whenever 
$m_d^{\circ i}(\lambda),m_d^{\circ j}(\lambda') \in \Lambda_k$ for some $i,j > 0$ 
then 
$m_d^{\circ i}(\lambda)=m_d^{\circ j}(\lambda')$.
(This is useful if we think of a dynamically preferred element in each $\Lambda_k$).
\bigskip
{\bf A.2.3 Linkage Relations.} 
We will say that two subsets $T$ and $T'$ of the circle $\bf R/Z$ are unlinked if 
they are contained in disjoint connected subsets of $\bf R/Z$, 
or equivalently, 
if $T'$ is contained in just one connected component of the complement ${\bf R/Z}-T$. 
(In particular $T$ and $T'$ must be disjoint.) 
If we identify $\bf R/Z$ with the boundary of the unit disk, 
an equivalent condition would be that the convex closures of these sets are pairwise disjoint. 
If $T$ and $T'$ are not unlinked then either $T \cap T'\ne \emptyset$ or there are elements $\theta_1,\theta_2 \in T$ and 
$\theta'_1,\theta'_2 \in T'$ such that the cyclic order can be written
$\theta_1, \theta'_1, \theta_2, \theta'_2,\theta_1$. 
In this second case we say that $T$ and $T'$ are {\it linked}.
More generally, a family ${\bf \Lambda}=\{\Lambda_1,\dots,\Lambda_n\}$ is an {\it unlinked family} if 
$\Lambda_1,\dots,\Lambda_n$ are pairwise unlinked.
Alternatively each $\Lambda_i$ is completely contained in a component of 
${\bf R/Z}-\Lambda_j$ for all $j \ne i$.
\medskip
The preceding definition has its motivation in the description of the dynamics of external rays in polynomial maps. 
Suppose the external rays $R_{\theta_i},R_{\psi_i}$ land at $z_i$ for $i=1,2$. 
If $z_1 \ne z_2$ then the sets $\{\theta_1,\psi_1\},\{\theta_2,\psi_2\}$ are unlinked, 
for otherwise the rays will cross each other.
The same argument applies if we consider rays supporting different Fatou components. 
But if we analyze linkage relations arising from rays supporting a Fatou component and rays that land at some point, 
we may get minor problems. 
Anyway, it is easy to see that even in this case the associated sets of arguments will be `almost' unlinked. 
(Compare condition (c.2) and as well as Proposition A.2.8 below.) 
\bigskip 
{\bf A.2.4 Weak linkage relations.} 
Consider families ${\bf {\cal F}}=\{{\cal F}_1,\dots,{\cal F}_n\}$ and ${\bf {\cal J}}=\{{\cal J}_1,\dots,{\cal J}_m\}$; 
we say that ${\bf \cal J}$ is {\it weakly unlinked to ${\bf \cal F}$ in the right} 
if we can chose arbitrarily small $\epsilon>0$ so that the family 
$\{{\cal F}_1,\dots,{\cal F}_n,\break 
{\cal J}_1-\epsilon,\dots,{\cal J}_m-\epsilon\}$ is unlinked. 
(Here ${\Lambda}-\epsilon=\{\lambda-\epsilon \hbox{ ({\it mod 1}) : }\lambda \in {\Lambda}\}$.)
In particular each family should be unlinked. 
Note that the definition allows empty families. 
To simplify notation we will simply say that 
{\it ``${\cal F}$ and ${\cal J}^-$ are unlinked"}. 
\bigskip
{\bf A.2.5 Formal Critical Portraits.}
Let ${\bf \cal F}=\{{\cal F}_1,\dots,{\cal F}_n\}$ and ${\bf \cal J}=\{{\cal J}_1,\dots,{\cal J}_m\}$ 
be two families of rational ({\it d-})prearguments. 
We say that the pair $({\bf {\cal F},{\cal J}})$ is a {\it degree d formal critical portrait} if the following conditions are satisfied. 
\smallskip
{\it 
(c.1) $d-1=\sum (\#({\bf \cal F}_k)-1)+\sum (\#({\bf \cal J}_l)-1)$
\smallskip
(c.2) ``${\cal F}$ and ${\cal J}^-$ are unlinked".
\smallskip
(c.3) Each family is hierarchic.
\smallskip
(c.4) Given $\gamma \in {\bf {\cal F}^\cup}$, there is an $i>0$ such that $m_d^{\circ i}(\gamma) \in {\bf {\cal F}^\cup_{per}}$.
\smallskip
(c.5) No $\theta \in {\bf {\cal J}^\cup}$ is periodic.}
\medskip
This set of conditions represent the simplest conditions satisfied 
by the critical marking of a postcritically finite polynomial. 
Condition (c.1) says that we have chosen the right number of arguments. 
Condition (c.2) means that the rays and extended rays determine sectors which do not cross each other, 
and that ${\bf \cal F}$ was constructed from arguments of left supporting rays. 
This reflects our decision to chose the supporting arguments as the rightmost possible argument of an external ray. 
Condition (c.3) reflects our choice of preferred rays. 
Condition (c.4) indicates that arguments in ${\bf \cal F}$ are related to Fatou critical points. 
Condition (c.5) indicates that arguments in ${\bf \cal J}$ are related to Julia set critical points.
Unfortunely there are formal critical portraits which do not correspond to a postcritically finite polynomial. 
In order to state necessary and sufficient conditions we need to study the 
dynamically defined partitions of the unit circle determined by these elements. 
\bigskip
{\bf A.2.6.} 
Given two families ${\bf {\cal F},{\cal J}}$ as above, 
we form a partition ${\bf \cal P}=\{L_1,\dots,L_d\}$ 
of the unit circle minus a finite number of points 
${\bf R/Z}-{\bf {\cal F}^\cup} - {\bf {\cal J}^\cup}$, in the following way. 
We consider two points $t,t' \in {\bf R/Z}-{\bf {\cal F}^\cup}-{\bf {\cal J}^\cup}$. 
By definition, $t,t'$ are {\it unlink equivalent} 
if they belong to the same connected component of ${\bf R/Z}- {\cal F}_i$ and ${\bf R/Z}- {\cal J}_j$, for all possible $i,j$.
Let $L_1,\dots,L_d$ be the resulting unlink equivalence classes with union ${\bf R/Z}-{\bf {\cal F}^\cup} - {\bf {\cal J}^\cup}$. 
It is easy to check that each $L_p$ is a finite union of open intervals with total length $1/d$.
\smallskip
Each element $L_i \in {\cal P}$ of the partition 
is a finite union $L_i=\cup (x_j,y_j)$ of open connected intervals.  
We define the sets $L_i^+=\cup [x_j,y_j)$ and\break 
$L_i^-=\cup (x_j,y_j]$. 
It is easy to see that both 
${\bf \cal P}^+=\{L_1^+,\dots,L_d^+\}$ and\break 
${\bf \cal P}^-=\{L_1^-,\dots,L_d^-\}$ 
are partitions of the unit circle. 
As every $\theta \in {\bf R/Z}$ belongs to exactly one set $L_k^+$, 
we define its {\it right address} $A^+(\theta)=L_k$. 
In an analogous way we define the {\it left address $A^-(\theta)$ of $\theta$}. 
We associate to every argument $\theta \in {\bf R/Z}$  
a {\it right symbol sequence} $S^+(\theta)=(A^+(\theta),A^+(m_d(\theta)),\dots)$ and 
a {\it left symbol sequence} 
$S^-(\theta)=(A^-(\theta),A^-(m_d(\theta)),\dots)$. 
Note that for all but a countable number of arguments $\theta \in {\bf R/Z}$ 
(namely the arguments present in the families and their iterated inverses), 
the left $S^-(\theta)$ and the right $S^+(\theta)$ symbol sequences agree. 
By $S(\theta)$ will be meant either (left or right) symbol sequence. 
\bigskip
{\bf A.2.7 Admissible Critical Portraits.}
Let ${\bf \cal F}=\{{\cal F}_1,\dots,{\cal F}_n\}$, ${\bf \cal J}=\{{\cal J}_1,\dots,{\cal J}_m\}$ 
be two families of rational ({\it d-})prearguments. 
We say that the pair $({\bf {\cal F},{\cal J}})$ is a 
{\it degree d admissible critical portrait} if 
$({\bf {\cal F},{\cal J}})$ is a degree $d$ formal critical portrait and the following two extra conditions are satisfied. 
\smallskip
{\it 
(c.6) Let $\gamma \in {\bf {\cal F}^\cup_{per}}$ and $\lambda \in {\bf R/Z}$, 
then $\lambda=\gamma$ if and only if $S^+(\gamma)=S^+(\lambda)$.
\smallskip
(c.7) Let $\theta \in {\cal J}_l$ and $\theta' \in {\cal J}_k$. 
If for some i, $S^-(m_d^{\circ i}(\theta))=S^-(\theta')$, then $m_d^{\circ i}(\theta) \in {\cal J}_k$.}
\bigskip
{\bf A.2.8 Proposition.}
{\it If $(P,{\cal F},{\cal J})$ is a critically marked polynomial, then $({\cal F},{\cal J})$ is an admissible critical portrait.} 
\medskip
Condition (c.6) indicates that arguments in ${\cal F}_l$ must support Fatou components. 
Condition (c.7) indicates that different elements in the family ${\cal J}$ are associated with different critical points. 
Now we can state the main result for critically marked polynomials as follows. 
\bigskip 
{\bf A.2.9 Theorem.}
{\it Let $({\bf {\cal F},{\cal J}})$ be a degree d admissible critical portrait. 
Then there is a unique monic centered postcritically finite polynomial P, 
with critical marking $(P,{\bf {\cal F},{\cal J}})$.}
\bigskip 
Now we should ask if conditions {\it (c.1)-(c.7)} represent a finite amount of information to be checked. 
This question is answered in a positive way by the following proposition. 
\medskip
{\bf A.2.10 Lemma.}
{\it Suppose $\theta$ and $\theta'$ have the same periodic (left or right) symbol sequence. 
Then $\theta$ and $\theta'$ are both periodic and of the same period.} 
\bigskip
{\bf A.2.11.} The next question that we ask is what kind of information 
about the Julia set can be gained by looking carefully into the combinatorics. 
For example, 
if can we determine if two rays land at the same point by only looking at their arguments. 
In fact, left symbol sequences 
contain all the information necessary to determine whether two rays land at the same point or not. 
This is done as follows. 
Suppose ${\cal J}_i=\{\theta_1,...,\theta_k\} \in {\cal J}$ with corresponding left symbol sequences 
$S^-(\theta_1),\dots,S^-(\theta_k)$. 
As we expect the rays with those arguments to land at the same critical point, 
we declare them ($i$-)equivalent; 
i.e, we write $S^-(\theta_\alpha) \equiv_i S^-(\theta_\beta)$. 
Then we set $\theta \approx \theta'$ either if 
$S^-(\theta)=S^-(\theta')$ or 
there is an $n \ge 0$ such that $A^-(m_d^{\circ j}(\theta))=A^-(m_d^{\circ j}(\theta'))$ for all $j < n$ and 
$S^-(m_d^{\circ n}(\theta)) \equiv_i S^-(m_d^{\circ n}(\theta'))$ for some $i$. 
This relation $\approx$ is not necessarily an equivalence relation, because transitivity may fail. 
To make this into an equivalence relation we say that $\theta \sim_l \theta'$ 
if and only if 
there are arguments $\lambda_0=\theta,\lambda_1,\dots,\lambda_m=\theta'$, 
such that $\lambda_0 \approx \dots \approx \lambda_m$. 
The importance of this equivalence relation is shown by the following proposition. 
\medskip
{\bf A.2.12 Proposition.} 
{\it Let $(P,{\bf {\cal F},{\cal J}})$ be a critically marked polynomial. 
Then $R_\theta$ and $R_{\theta'}$ land at the same point if and only if $\theta \sim_l \theta'$.} 
\bigskip
{\bf A.2.13 Corollary.} 
{\it 
The symbol sequence $S^-(\theta)$ is a periodic sequence of period $m$ 
if and only if 
the landing point of the ray $R_\theta$ has period $m$.} 
\endofproof
\bigskip
\inchap
\centerline{\bigfont Appendix B}
\centerline{\bigfont Finite Cyclic Expanding Maps.}
\inchap
\insec
\centerline{\medfont 1. Expanding Maps.}
\insec
{We consider a finite cyclic set $X$, 
and a degree $n \ge 2$ orientation preserving map  
$f:X \mapsto f(X) \subset X$. 
We will study under which conditions we can assign an {\it argument} $\phi(p)$ to every point $p \in X$ 
such that the induced map becomes multiplication by $n$.}  
\bigskip
{\bf 1.1.} 
Let $k \ge 1$ and $n \ge 2$. 
Consider a finite cyclicly ordered set $X=\{p_1,\dots,p_{kn}\}$ with $kn$ elements. 
The {\it cyclic order} can be realized as a {\it successor} 
function $Suc_X(p_i)=p_{i+1}$ with the convention $p_{kn}=p_0$. 
Given $Y \subset X$ there is an induced order in $Y$, and therefore a successor function $Suc_Y:Y \mapsto Y$. 
We consider a degree $n \ge 2$ orientation preserving map $f:X \mapsto f(X) \subset X$. 
By this we mean a function $f$ with the property that $f(p_i)=f(p_j)$ if and only if $i \equiv j$ {\it (mod k)}; 
and such that $f(Suc_X(p))=Suc_{f(X)}(f(p))$. 
It follows that $f$ is an $n^{th}$-fold cover of its image. 
Note that because $f$ is a degree $n$ cover and order preserving, 
for every $p \in X$, 
the restriction of $f$ to the set $\{p,Suc_X(p),\dots,Suc_X^{\circ k-1}(p)\}$ is one to one and onto $f(X)$. 
\medskip
Given a cyclicly ordered set $X$ as above, 
we define the {\it ordered distance $d_X(p_1,p_2)$ between two points $p_1,p_2 \in X$}, 
as the minimal $m$ for which $p_2=Suc^{\circ m}(p_1)$. 
Thus, the ordered distance between two points is always less than $kn$.
It follows easily that $f(p_1)=f(p_2)$ if and only if $d_X(p_1,p_2)$ is a multiple of $k$. 
\bigskip
Given three points $p_1,p_2,p_3$ and numbers $0 \le m \le m'< kn$, 
with $m=d_X(p_1,p_2)$, $m'=d_X(p_1,p_3)$, 
we write $p_1 \le p_2 \le p_3$. 
If in addition $m <m'$ we write $p_1 \le p_2 < p_3$. 
\medskip
{\bf 1.2 Lemma.}
{\it Suppose $p_1 \le p_2 \le p_3 < p_1$. 
Then $d_X(p_1,p_3)\break 
= d_X(p_1,p_2)+ d_X(p_2,p_3)$.} 
\medskip
{\bf Proof.} 
Completely trivial. 
\endofproof
\bigskip
{\bf 1.3 Remark.} 
Even if we are considering two orders (one in $X$ and that induced in $f(X)$), 
we will only be considering the ordered distance of $X$. 
In other words if  
$p_1,p_2 \in f(X)$, 
the ordered distance $d_X(p_1,p_2)$ is always measured in $X$. 
\bigskip
{\bf 1.4 Definition.} 
We say that $f:X \to X$ as above is {\it expanding}, if 
given $p_1,p_2$ periodic 
\smallskip
{\it ($\star$) there exists $l \ge 0$ such that 
$d_X(f^{\circ l}(p_1),f^{\circ l}(p_2)) \ne 1$. 
} 
\medskip
In other words, if two periodic points are consecutive, 
the distance between them eventually increases.  
From the facts that $d_X(p_1,p_2)<k$ implies $f(p_1) \ne f(p_2)$, 
and every point is eventually periodic, 
we can easily deduce that for an expanding map, 
condition {\it ($\star$)} is also satisfied for every pair of different points. 
\bigskip
{\bf 1.5} 
Given a finite cyclicly ordered set $X$ and a degree $n \ge 2$ orientation preserving map $f$, 
we say that $f:X \to X$ can be {\it angled} 
if there is an order preserving embedding $\phi:X \mapsto {\bf R/Z}$, 
such that $n\phi(p) \equiv \phi(f(p)) \hbox{\it (mod 1)}$. 
Of course, an angled function is expanding. 
\bigskip
{\bf 1.6 Remark.} 
If we reverse the order in all the definitions above 
(i.e, if we replace the successor function by a predecessor function $Pre_X$), 
all the definitions above make sense. 
In particular if $\phi_S$,$\phi_P$ are the angle functions for these two orders then clearly 
$\phi_S+\phi_P \equiv 1$. 
\bigskip
{\bf 1.7 Proposition.} 
{\it Let X be a finite cyclic set and f an orientation preserving degree $n\ge 2$ map. 
Then $f:X \to X$ is angled if and only if is expanding.}
\medskip
{\bf Proof.}
Being angled implies being expanding as remarked above. 
We prove the converse in several steps. 
\smallskip 
{\it Step 1:} We can assume without loss of generality that $f$ has a fixed point. 
In fact, if there is no fixed point, then  
$f(X)$ has at least two elements. 
We define a function 
$g:X \to \{1,\dots,kn-1\}$ by the formula $g(x)=d_X(x,f(x))$. 
It follows easily from Lemma 1.2 that whenever $i \equiv j$ {\it (mod k)} then 
$g(x_i) \equiv g(x_j) + d_X(x_i,x_j)$ {\it (mod kn)}. 
Therefore for $y \in f(X)$, there is a unique 
$x_i \in f^{-1}(y)$ for which 
$k(n-1)<g(x_i)<kn$ 
(in fact, $g(x_i)=k(n-1)$ would imply that $x_{i+k}$ is a fixed point). 
Let $d$ be the maximum of $g$. 
Among all $x$ with $g(x)=d$ take one for which $g(Suc_X(x))<d$. 
It follows easily that the cyclic order can be written as 
$$
f(x) < x < Suc_X(x) < f(Suc_X(x))=Suc_{f(X)}(f(x)) < f(x). 
$$
To simplify notation, we rewrite $X$ as $\{p_0=x,p_1,\dots,p_{kn-1}\}$. 
We insert a new point $q_i$ between every pair $p_{ki}$ and $p_{ki+1}$. 
All of this new points will be mapped to $q_0$.  
In this way, 
we have a degree $n \ge 2$ orientation preserving map which is an extension of the original one. 
\smallskip
We must verify that this map is expanding. 
The only new periodic point included is $q_0$. 
The expanding property is obviously verifies if $Suc_X(q_0)$ is periodic: 
if $d_X(f(q_0),f(Suc_X(q_0)))=1$ then $Suc_X(q_0)$ is a fixed point, in contradiction to what was assumed. 
If $Pre_X(q_0)$ is periodic the result follows analogously. 
\medskip
{\it Step 2:}
We assign an argument to each point in $X$ as follows. 
Let $q_0 < q_1 < q_{n-1} < q_0$ be all points which map to the fixed point $q_0$. 
We assign to $q_i$ the argument $i/n$ for $i=0,\dots,n-1$. 
For an arbitrary point $x \in X$, 
we dynamically find its numerical expansion in base $n$.
\medskip
{\it Step 3: The assignment is order preserving.} 
Because the function is $n$ to one order preserving, we may introduce inverse iterates of the fixed point. 
Thus, we may assume that given $m$ there are in the cyclic order different values $\{q_0=q,\dots,q_{m^n-1}\}$ with the property that $f^{\circ m}(q_i)=q_0$. 
Taking $m$ big enough the result follows.
\medskip
{\it Step 4: Different points are assigned different arguments.} 
Consider a set $\{x_1,\dots,x_l\}$ of maximal cardinality to which equal periodic 
base $n$ expansion is associated. 
Clearly all $x_i$ are periodic. 
Furthermore, if $l>1$ we have for all $m \ge 0$ 
$d_X(f^{\circ m}(x_1),f^{\circ m}(x_2))=1$ because of maximality. 
But this contradicts the expanding condition. 
There is a case in which this argument does not apply. 
Suppose that in applying step 2, there is an argument to which the decimal expansion 
$0.n-1,n-1,\dots$ is assigned.  
In this\break 
case we reverse the order, and apply the same argument to derive a\break 
contradiction.  
\endofproof
\insec
\centerline{\medfont 2. Finding the Coordinates.}
\insec
{ 
Consider an integer $n>1$, 
and denote by $m_n$ multiplication by $n$ modulo 1. 
From the dynamical point of view the election of $0$ as the origin is arbitrary in the sense that any dynamically property present at a point $x \in {\bf R/Q}$, 
is also present at $x+j/(n-1)$. 
In this way, with the knowledge of the dynamical behavior of a point $x$, 
the natural question is not which is the value of $x$, but that of $m_{n-1}(x)$.}
\bigskip
{\bf 2.1}
For $n>1$ define 
$\delta_n:{\bf R/Z} \to \{0,\dots,n-1\}$ by 
$$
\delta_n(x)=i \hskip 0.2in \hbox{ if \hskip 0.2in } 
m_{n-1}(x) \in [{i \over n},{i+1 \over n}).
$$
In other words, if we take $m_{n-1}(x)$, 
we define $\delta_n(x)$ as the integer part of $n(m_{n-1}(x))$. 
\bigskip
{\bf 2.2 Remark.} 
It follows that 
$\delta_n(x)$ is the number of inverses of $m_n(x)$ (other than $x$) in the cyclicly counterclockwise oriented interval $(x,m_d(x))$ 
(if $x$ is fixed this interval is interpreted to be empty).
To see this, we rewrite $m_{n-1}(x)$ as $m_{n}(x)-x$ {\it (mod n)}. 
In this way, 
$\delta_n(x)$ counts the number of intervals of size $1/n$ to be found in $(x,m_d(x))$. 
The claim follows easily. 
\bigskip
{\bf 2.3 Example.} 
Consider with $n=3$ the point $x=1/5$. 
We have $m_2(x)=2/5$ and $1 \le 3(2/5) <2$; 
so by definition $\delta_3(1/5)=1$. 
Note also that $\delta_3(1/5+1/2)=1$, which is not a surprise because 
$m_2(x)=m_2(x+1/2)$ for all $x \in {\bf R/Z}$. 
\bigskip
{\bf 2.4 Lemma.} 
{\it Let $n>1$ and $x \in {\bf R/Z}$, then
$$
m_{n-1}(x)={1 \over n}\sum_{i=0}^\infty {\delta_n(m_n^{\circ i}(x)) \over n^i}.
$$
}
\medskip
{\bf Proof.} 
We successively subdivide the interval $I_j=[{j \over {n-1}},{{j+1} \over {n-1}})$ in $n$ semiopen intervals. 
This determines a parametrization of the interval $I_j$ by symbol sequences in the symbol space $\{0,\dots,n-1\}$ 
(not allowing any symbol sequence with tail $(n-1,n-1,\dots)$). 
A point $x \in I_j$ has symbol sequence $S_0,S_1,\dots$ if and only if 
$x={j \over {n-1}}+ {1 \over (n-1)n} \sum_{i=0} {S_i \over n^i}$. 
Therefore, $m_{n-1}(x)= {1 \over n}\sum_{i=0}^\infty {S_i \over n^i}$ 
(and all reference to the initial interval $I_j$ is lost). 
The result follows as $S_i=\delta_n(m_n^{\circ i}(x))$ by construction.
\endofproof
\bigskip 
{\bf 2.5 Coordinates for expanding maps.} 
We return to the case described in $\S$1. 
It follows by definition of covering map and Remark 2.2 that 
$\delta_n(x)$ equals the integer part of 
${d_X(x,f(x))}\over {k}$. 
Thus, according to Lemma 2.4, 
$m_{n-1}(x)$ is independent of the coordinate assigned in Proposition 1.7. 
Furthermore, we have proved the following. 
\bigskip 
{\bf 2.6 Theorem.} 
{\it Let X be a finite cyclic ordered set and $f:X \mapsto X$ be a degree n orientation preserving expanding map. 
Then f can angled in exactly $n-1$ ways.}
\bigskip
{\bf 2.7 Example.} 
(Compare Figure B.1.) 
Let $X$ be the cyclic set shown in Figure $B.1$. 
(The notation is justified by the dynamics.) 
We consider a map $f:X \to X$ for which 

\centerline{$f(A)=f(A')=f(A'')=B$} 
\centerline{$f(B)=f(B')=f(B'')=C$} 
\centerline{$f(C)=f(C')=f(C'')=D$} 
\centerline{$f(D)=f(D')=f(D'')=A$}. 

The unique periodic orbit is given by $A \mapsto B \mapsto C \mapsto D \mapsto A$. 
This map is clearly expanding. 
According to Remark 2.5, we have that   
$\delta(A)=1$, $\delta(B)=0$, $\delta(C)=1$, and $\delta(D)=2$.
Using Lemma 2.4, we can easily find the base 3 expansions of $m_2(A)=0.\overline{1012}$. 
It follows that $m_2(A)=2/5$, and therefore A takes value either $1/5$ or $7/10$.
\medskip
\centerline{
\psfig{figure=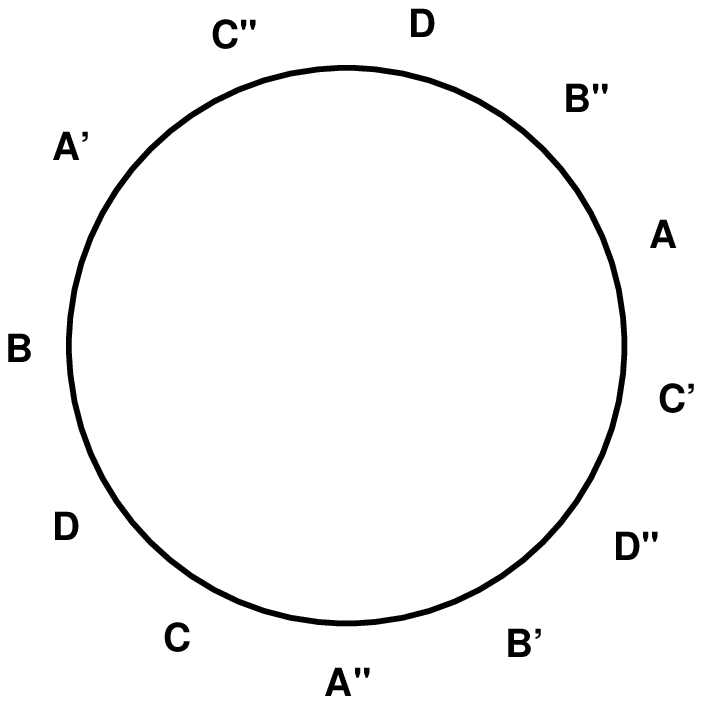,width=2.5in}
}
\medskip
\centerline{\sl Figure B.1}
\bigskip 
{\bf 2.8 Corollary.}
{\it Let x be periodic under $m_n$, and denote by ${\cal O}(x)$ its orbit. 
Then $m_{n-1}(x)$ is uniquely determined by the cyclic order of $m_n^{-1}({\cal O}(x))$.}
\medskip
{\bf Proof.} 
This follows directly from Remark 2.2 and Lemma 2.4. 
(Compare also Example 2.7). 
\endofproof 
\vfill\eject 
\inchap
\centerline{\bigfont References.}
\inchap
[BFH] B. Bielefeld, Y.Fisher, J. Hubbard, The Classification of Critically Preperiodic Polynomials as Dynamical Systems; 
Journal AMS 5(1992)pp. 721-762.
\smallskip
[DH1] A.Douady and J.Hubbard, \'Etude dynamique des polyn\^omes complexes, part I; 
Publ Math. Orsay 1984-1985. 
\smallskip
[DH2] A.Douady and J.Hubbard, A proof of Thurston's Topological Characterization of Rational Maps; Preprint, Institute Mittag-Leffler 1984. 
\smallskip
[F] Y.Fisher, Thesis; Cornell University, 1989. 
\smallskip
[GM] L.Goldberg and J.Milnor, Fixed Point Portraits; 
Ann. scient. \'Ec. Norm. Sup., $4^e$ s\'erie, t. 26, 1993, pp 51-98 
\smallskip
[L] P.Lavaurs, These; Universit\'e de Paris-Sud Centre D'Orsay; 1989. 
\smallskip
[M] J.Milnor, Dynamics in one complex variable: Introductory Lectures; 
Preprint \#1990/5 IMS SUNY@StonyBrook.
\smallskip
[P1] A.Poirier, On Postcritically Finite Polynomials, Part One: Critical Portraits; 
Preprint \#1993/5 IMS SUNY@StonyBrook.
\smallskip
[P2] A.Poirier, Thesis; On Postcritically Finite Polynomials;\break  
SUNY@StonyBrook 1993.
\smallskip

\end